\documentclass[final]{siamltex}

\usepackage{times}
\usepackage{epsfig}
\usepackage{graphicx}
\usepackage{graphics} 
\usepackage{wrapfig}
\usepackage{amsmath}
\usepackage{amssymb}
\usepackage{dsfont}
\usepackage{xspace}
\usepackage{color}

\newcommand{\object}{{\mathcal{O}}}
\newcommand{\refobject}{{\mathcal{O}_{\mbox{\tiny ref}}}}

\newcommand{\curve}{\object}
\newcommand{\V}{{\mathcal{V}}}

\newcommand{\ikm}{{k-1}}
\newcommand{\ik}{{k}}
\newcommand{\compDom}{{D}}

\newcommand{\manifold}{{\mathcal{M}}}
\newcommand{\setof}[2]{{\{{#1}\,|\,{#2}\}}}

\newcommand{\metric}{{\mathcal{G}}}
\newcommand{\pathlength}{{\mathrm{L}}}
\newcommand{\pathenergy}{{\mathrm{E}}}
\DeclareMathOperator{\argmin}{argmin}
\newcommand{\Log}[1]{{(\textstyle\frac1{#1}\mathrm{LOG})}}
\newcommand{\Exp}[1]{{\mathrm{EXP}^{#1}}}
\newcommand{\ExpO}[2]{{\mathrm{EXP}_{#2}^{#1}}}

\newcommand{\W}{{\mathcal{W}}}
\newcommand{\dissdensity}{\mathbf{diss}}
\newcommand{\F}{{\mathcal{F}}}

\newcommand{\R}{{\mathds{R}}}

\renewcommand{\d}{{\mathrm{d}}}
\newcommand{\D}{{\mathcal{D}}}

\newcommand{\dist}{{\mathrm{dist}}}
\newcommand{\tr}{{\mathrm{tr}}}
\newcommand{\id}{{\mathrm{id}}}
\newcommand{\Id}{{\mathds{1}}}

\newcommand{\notinclude}[1]{}

\makeatletter
\DeclareRobustCommand\onedot{\futurelet\@let@token\@onedot}
\def\@onedot{\ifx\@let@token.\else.\null\fi\xspace}

\def\eg{\emph{e.g}\onedot} 
\def\ie{\emph{i.e}\onedot}

\def\etal{\emph{et al}\onedot}
\def\conf{\emph{cf}\onedot}
\makeatother

\begin{document}

\title{Discrete geodesic calculus in the space of viscous fluidic objects}


\author{Martin Rumpf\thanks{Bonn University, Endenicher Allee 60, D-53115 Bonn, Germany ({\tt Martin.Rumpf@ins.uni-bonn.de}).}
\and Benedikt Wirth\thanks{Courant Institute, New York University, 251 Mercer Street, New York, NY 10012, USA ({\tt Benedikt.Wirth@cims.nyu.edu}).}}

\maketitle

\begin{abstract}
Based on a local approximation of the Riemannian distance on a manifold by a computationally cheap dissimilarity measure, a time discrete geodesic calculus is developed, and applications to shape space are explored.
The dissimilarity measure is derived from a deformation energy whose Hessian reproduces the underlying Riemannian metric, and it is used to define length and energy of discrete paths in shape space.
The notion of discrete geodesics defined as energy minimizing paths gives rise to
a discrete logarithmic map, a variational definition of a discrete exponential map, and a time discrete parallel transport.
This new concept is applied to a shape space in which shapes are considered as boundary contours of
physical objects consisting of viscous material. The flexibility and computational efficiency of the approach is
demonstrated for topology preserving shape morphing,
the representation of paths in shape space via local shape variations
as path generators, shape extrapolation via discrete geodesic flow, and the transfer of geometric features.
\end{abstract}



\pagestyle{myheadings}
\thispagestyle{plain}
\markboth{MARTIN RUMPF AND BENEDIKT WIRTH}{DISCRETE GEODESIC CALCULUS IN THE SPACE OF VISCOUS FLUIDIC OBJECTS}

\section{Introduction}
Geodesic paths in shape space allow to define smooth and in some sense geometrically or physically natural
connecting paths $\object(t)$, $t\in[0,1]$, between two given shapes $\object(0),\object(1)$, or they enable the extrapolation of a path from an initial shape $\object(0)$ and an initial shape variation $\delta\object$ which encodes the path direction.
Applications include shape modeling in computer vision \cite{KlSrMi04,KiMiPo07},
computational anatomy, where the morphing path establishes correspondences between a patient and a template \cite{BeMiTrYo02,MiTrYo03},
shape clustering based on Riemannian distances \cite{SrKlJo11},
as well as shape statistics \cite{FlLuPi04,FuSc08}, where geodesic paths in shape space transport information from the observed shapes into a common reference frame in which statistics can be performed.

As locally length minimizing paths, geodesic paths require to endow the space of shapes with a Riemannian metric
which encodes the preferred shape variations.
There is a rich diversity of Riemannian shape spaces in the literature.
Kilian\,\etal compute isometry invariant geodesics between consistently triangulated surfaces \cite{KiMiPo07},
where the Riemannian metric measures stretching of triangle edges,
while the metric by Liu\,\etal also takes into account directional changes of edges \cite{LiShDi10}.
\begin{figure}
\setlength{\unitlength}{.145\linewidth}
\begin{picture}(7,1.7)
\put(0,1){\includegraphics[width=.95\unitlength]{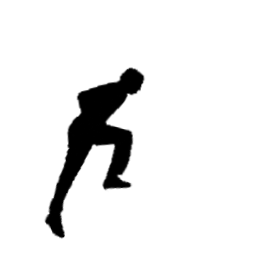}}
\put(1,1){\includegraphics[width=.95\unitlength]{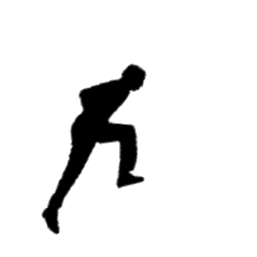}}
\put(0.72,1.4){\Large \color{red} $\stackrel{\zeta_1}{\longrightarrow}$}
\put(2,1){\includegraphics[width=.95\unitlength]{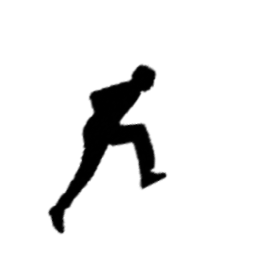}}
\put(3,1){\includegraphics[width=.95\unitlength]{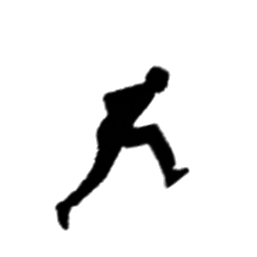}}
\put(4,1){\includegraphics[width=.95\unitlength]{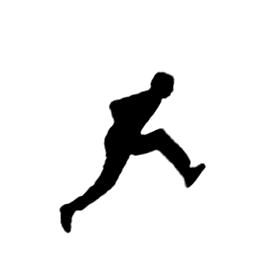}}
\put(5,1){\includegraphics[width=.95\unitlength]{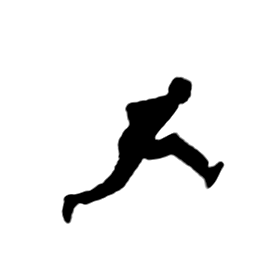}}
\put(6,1){\includegraphics[width=.95\unitlength]{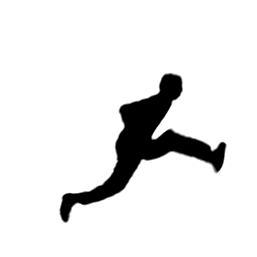}}
\put(0,0){\includegraphics[width=.95\unitlength]{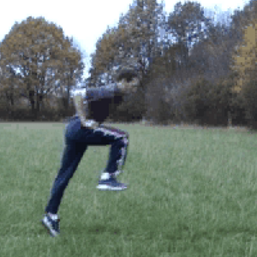}}
\put(1,0){\includegraphics[width=.95\unitlength]{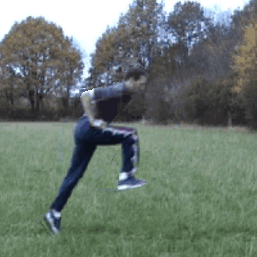}}
\put(2,0){\includegraphics[width=.95\unitlength]{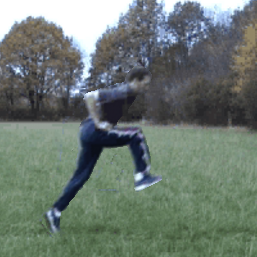}}
\put(3,0){\includegraphics[width=.95\unitlength]{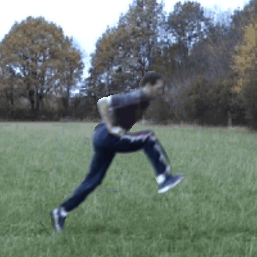}}
\put(4,0){\includegraphics[width=.95\unitlength]{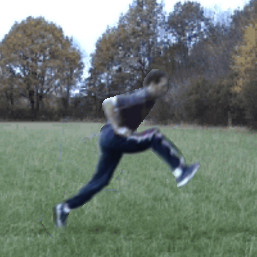}}
\put(5,0){\includegraphics[width=.95\unitlength]{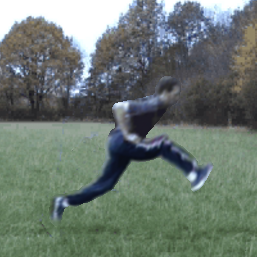}}
\put(6,0){\includegraphics[width=.95\unitlength]{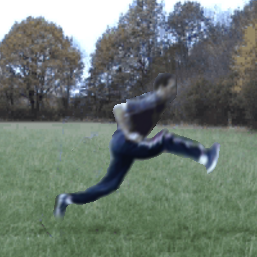}}
\end{picture} \label{fig:teaser}
\caption{Top: Given the first shape on the left and an initial variation $\zeta_1$ described by the second shape,
a discrete geodesic path is extrapolated.
Bottom: The texture of a video frame can be transported along with the resulting geodesic flow.}
\end{figure}

For planar curves, different Riemannian metrics have been devised, including the $L^2$-metric on direction and curvature functions \cite{KlSrMi04},
the $L^2$-metric on stretching and bending variations \cite{SrJaJo06},
as well as curvature-weighted $L^2$- or Sobolev-type metrics \cite{MiMu04,SuYeMe07},
some of which allow closed-form geodesics \cite{YoMiSh08,SuMeSo11}.
A variational approach to the computation of geodesics in the space of planar Jordan curves has been proposed by Schmidt \etal in \cite{ScClCr06}. The extrapolation of geodesics in the space of curves incorporating translational, rotational, and scale invariance has been investigated by Mennucci \etal \cite{MeSo10}.
A Riemannian space of non-planar elastic curves has very recently been proposed by Srivastava \etal \cite{SrKlJo11}.

In the above approaches, the shape space is often identified
with a so-called pre-shape space of curve parameterizations over a 1D domain (or special representations thereof)
modulo the action of the reparameterization group.
It is essential that the metric on the pre-shape space is invariant under reparameterization
or equivalently that reparameterization represents an isometry in the pre-shape space
so that the Riemannian metric can be inherited by the shape space.
Such reparameterization-invariant metrics can also be defined on the space of parameterized 2D surfaces \cite{BaBr11,KuKlGo11}.
For certain representations of the parameterization one is lead to a very simple form of the metric, \eg an $L^2$-type metric \cite{KuKlDi10}.

The issue of reparameterization invariance does not occur when the mathematical description of the shape space is not based on parameterizations,
which often simplifies the analysis (and is also the approach taken here).
When warping objects in $\R^d$, a shape tube in $\R^{d+1}$ is formed.
Zol\'esio investigates geodesic in terms of shortest shape tubes \cite{Zo03}.
The space of sufficiently smooth domains $\object\subset\R^d$ can be assigned a Riemannian metric
by identifying the tangent space at $\object$ with velocity fields $v:\object\to\R^d$ and defining a metric on these.
Dupuis\,\etal employ a metric 
$$\metric(v,v)=\int_\compDom L v \cdot v \,\d x$$
for a higher order elliptic operator $L$ on some computational domain $\compDom\subset\R^d$ \cite{DuGrMi98},
ensuring a diffeomorphism property of geodesic paths.
A corresponding geodesic shooting method has been implemented in \cite{BeMiTr05}.
Fuchs\,\etal propose a viscous-fluid based Riemannian metric \cite{FuJuScYa09}.
Fletcher and Whitaker employ a similar metric on pullbacks of velocity fields onto a reference shape \cite{FlWh06}.
Miller and Younes consider the space of registered images as the product space of the Lie group of diffeomorphisms and image
maps. They define a Riemannian metric using sufficiently regular elliptic operators on the diffeomorphism-generating velocity fields,
which may also depend on the current image \cite{MiYo01}.
A morphing approach based on the concept of optimal mass transport has been proposed by Haker \etal  \cite{HaTaKi01,ZhYaHa07}.
An image or a shape is viewed as mass density,
and for two such densities $\rho_0,\rho_1:\compDom\to\R$ the Monge--Kantorovich functional 
$$\int_\compDom|\psi(x)-x|^2\rho_0(x)\,\d x$$ is minimized over all mass preserving mappings $\psi:\compDom\to\compDom$, \ie mappings with $\rho_0=\rho_1\circ\psi\det\nabla\psi$.
A morphing path then is given by $\rho(t)=\rho_0\circ\psi(t)^{-1}\det \nabla \psi(t)^{-1}$ for $\psi(t)=t\psi+(1-t)\id$, $t\in[0,1]$.
Like for our approach there is a continuum-mechanical interpretation of minimizing the action of an incompressible fluid flow \cite{BeBr00}, however, the flow typically does neither preserve local shape features or isometries nor the shape topology.

Very often, geodesics in shape space are approached via the underlying geodesic evolution equation,
and geodesics between two shapes are computed by solving this ODE within a shooting method framework \cite{KlSrMi04,BeMiTr05,BaBr11}.
An alternative approach exploits the energy-minimi\-zing property of geodesics:
Schmidt \etal perform a Gau\ss-Seidel type fixed-point iteration which can be interpreted as a gradient descent on the path energy,
and Srivastava \etal derive the equations of a gradient flow for the path energy which they then discretize \cite{SrKlJo11}.
In contrast, we employ an inherently variational formulation where geodesics are defined as minimizers of a time discrete path energy. Discrete geodesics are then defined consistently as minimizers of a corresponding discrete energy.

In this paper we start from this time discretization and consistently develop a time discrete geodesic calculus in shape space. 
The resulting variational discretization of the basic Riemannian calculus consists of an exponential map, a logarithmic map, parallel transport, and finally an underlying discrete connection. 
To this end, we replace the exact, computationally expensive Riemannian distance by a relatively cheap but consistent dissimilarity measure.
Our choice of the dissimilarity measure not only ensures consistency for vanishing time step size but also
a good representation of shape space geometry already for coarse time steps. For example, rigid body motion invariance is naturally incorporated in this approach. We illustrate this approach on a shape space consisting of homeomorphic
viscous-fluid objects and a corresponding deformation-based dissimilarity measure.

Different from most approaches, which first discretize in space
via the choice of a parameterization, a set of control points, or a mesh, and then solve the resulting transport equations
by suitable solvers for ordinary differential equations (see the discussion above), our time discretization is defined
on the usually infinite dimensional shape space. It results from a consistent transfer
of time continuous to time discrete variational principles. Thereby, it leads 
to a collection of variational problems on the shape space, 
which in our concrete implementation of the proposed calculus consists of non-parameterized, volumetric objects.

Let us also already mention a further remarkable conceptual difference.
The way the time discrete geodesic calculus is introduced differs substantially from the way
the time continuous counterpart is usually developed.
In classical Riemannian differential geometry one first defines a connection $(v,w) \mapsto \nabla_v w$
for two vector fields $v$ and $w$ on a manifold $\manifold$. With the connection at hand a tangent vector $w$ can be transported parallel along a path with motion field $v$ solving $\nabla_v w =0$. Studying those paths where the motion field itself is transported parallel along the path (\ie it solves the ODE $\nabla_v v =0$) one is led to geodesics.
Next, the exponential map is introduced via the solution of the above ODE for varying initial velocity. 
Finally, the logarithm is obtained as the (local) inverse of the exponential map. 

In the time discrete calculus we start with a time discrete formulation of path length and energy and then define discrete
geodesics as minimizers of the discrete energy. Evaluating the initial step of a discrete geodesic path as the discrete counterpart of the initial velocity we are led to the discrete logarithm. Then, the discrete exponential map is defined as the inverse of the discrete logarithm. Next, discrete logarithm and discrete exponential allow to define a discrete parallel
transport based on the construction of a sequence of approximate Riemannian parallelograms (commonly known as Schild's ladder \cite{EhPiSc72}). Finally, with the discrete parallel transport at hand, a time discrete connection can be defined.

Let us note that the approximation of parallel transport in shape space via Schild's ladder
has also been used in the context of the earlier mentioned flow of diffeomorphism approach \cite{PeLo11,LoAyPe11}.
In our discrete framework, however, the notion of discrete parallel transport is directly derived from the parallelogram construction,
consistently with the overall discrete approach to geodesics.

A related approach for time discrete geodesics has been presented in an earlier paper by Wirth \etal \cite{WiBaRu10}.
In contrast to \cite{WiBaRu10}, we here do not restrict ourselves to the computation of geodesic paths between two shapes but
devise a full-fledged theory of discrete geodesic calculus (\conf Figure~\ref{fig:teaser}). Furthermore, different from that approach we ensure topological consistency and describe shapes solely via deformations of reference objects instead of treating deformations and level set representations of shapes simultaneously as degrees of freedom,
which in turn strongly simplifies the minimization procedure.


The paper is organized as follows. In Section \ref{sec:viscous} we introduce a special model for a shape space, 
the space of viscous fluidic objects, to which we restrict our exposition of the geodesic calculus. Here, in the light of the discrete shape calculus to be developed, we will review the notion of discrete path length and discrete path energy. After these preliminaries the actual time discrete calculus consisting of a discrete logarithm, a discrete exponential and a discrete parallel transport together with a discrete connection is introduced and discussed in Section \ref{sec:timediscrete}.
Then, Section \ref{sec:numerics} is devoted to the numerical discretization via characteristic functions and a parameterization via deformations over reference paths.  
Finally, we draw conclusions in Section \ref{sec:conclusions}.

\section{A space of volumetric objects and an elastic dissimilarity measure}
\label{sec:viscous}
To keep the exposition focused we restrict ourselves to a specific shape model, where shapes are represented by
volumetric objects which behave physically like viscous fluids.
In fact, the scope of the variational discrete geodesic calculus extends beyond this concrete shape model.
We refer to Section\,\ref{sec:conclusions} for remarks on the application to more general shape spaces.

\subsection{The space of viscous-fluid objects}
\label{sec:viscousshapespace}
\begin{figure}
\setlength{\unitlength}{.125\linewidth}
\begin{picture}(8,2)(.15,-.05)
\put(0,1){\includegraphics[width=1.2\unitlength]{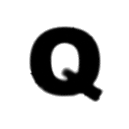}}
\put(1,1){\includegraphics[width=1.2\unitlength,angle=-9,origin=c]{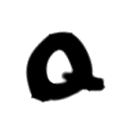}}
\put(2,1){\includegraphics[width=1.2\unitlength,angle=7,origin=c]{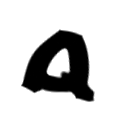}}
\put(3,1){\includegraphics[width=1.2\unitlength,angle=-0,origin=c]{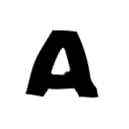}}
\put(4,1){\includegraphics[width=1.2\unitlength]{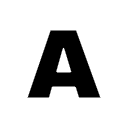}}
\put(4.8,1){\includegraphics[width=1.2\unitlength,angle=32,origin=c]{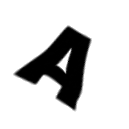}}
\put(6,1){\includegraphics[width=1.2\unitlength,angle=-11,origin=c]{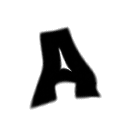}}
\put(7,1){\includegraphics[width=1.2\unitlength,angle=22,origin=c]{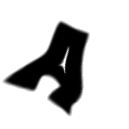}}
\put(0,.1){\includegraphics[width=1.2\unitlength]{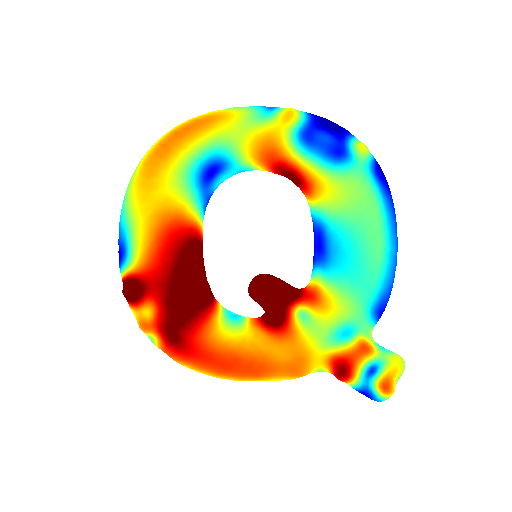}}
\put(1,.1){\includegraphics[width=1.2\unitlength,angle=-9,origin=c]{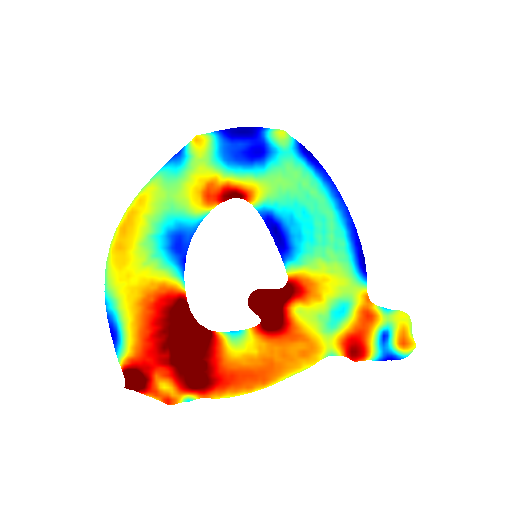}}
\put(2,.1){\includegraphics[width=1.2\unitlength,angle=7,origin=c]{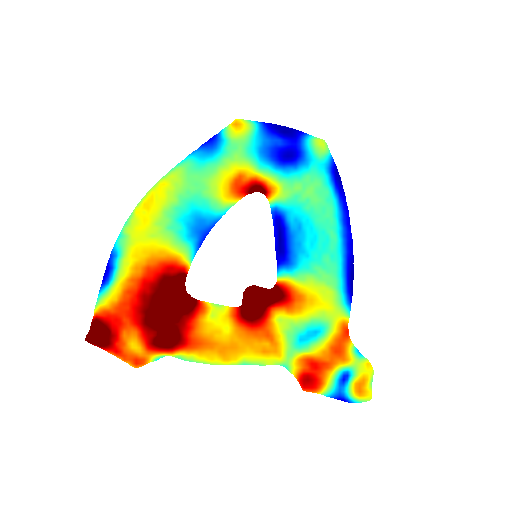}}
\put(3,.1){\includegraphics[width=1.2\unitlength,angle=-0,origin=c]{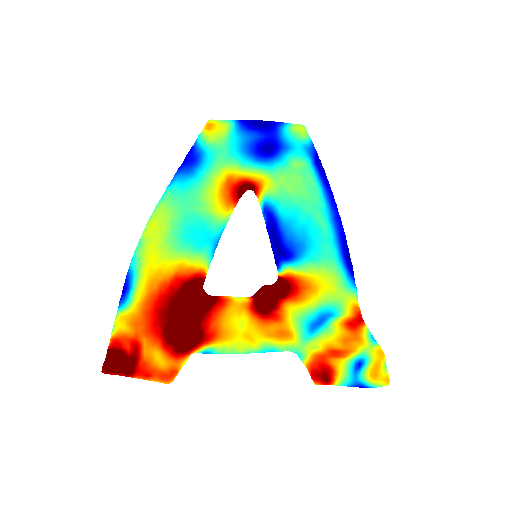}}
\put(4,.1){\includegraphics[width=1.2\unitlength]{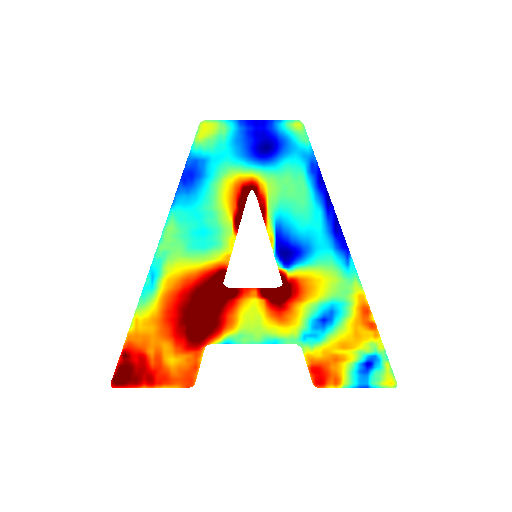}}
\put(4.8,.1){\includegraphics[width=1.2\unitlength,angle=32,origin=c]{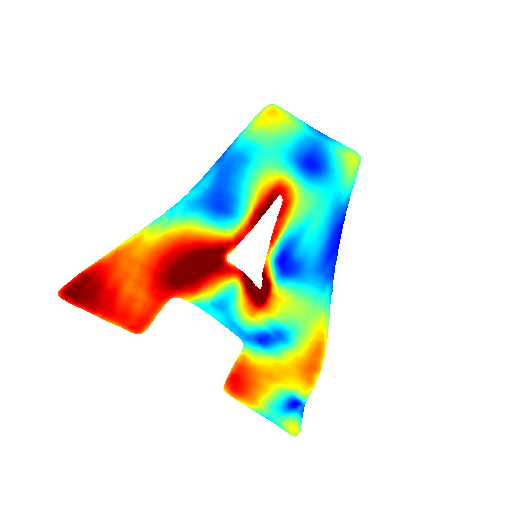}}
\put(6,.1){\includegraphics[width=1.2\unitlength,angle=-11,origin=c]{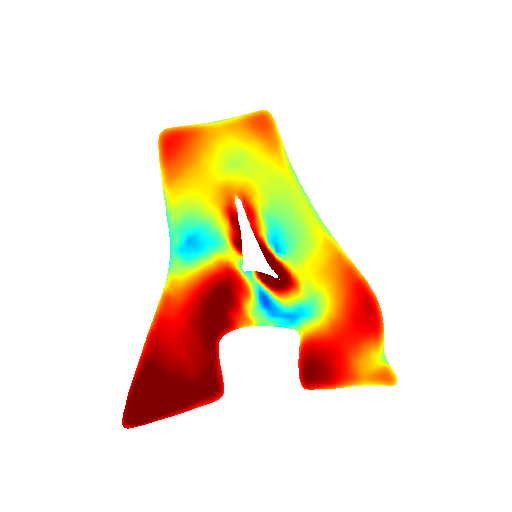}}
\put(7,.1){\includegraphics[width=1.2\unitlength,angle=22,origin=c]{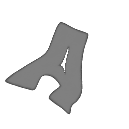}}
\put(5.1,-.1){\line(0,1){2.05}}
\put(1.8,0.05){\small geodesic matching}
\put(5.3,0.05){\small extrapolation}
\put(1.3,.1){\vector(-1,0){1}}
\put(3.85,.1){\vector(1,0){1}}
\put(7,.1){\vector(1,0){1}}
\end{picture}
\caption{Discrete geodesic between the letters Q and A, extrapolated beyond A.
Colors indicate the local rate of dissipation (from blue, low, to red, high).}
\label{fig:QA}
\end{figure}
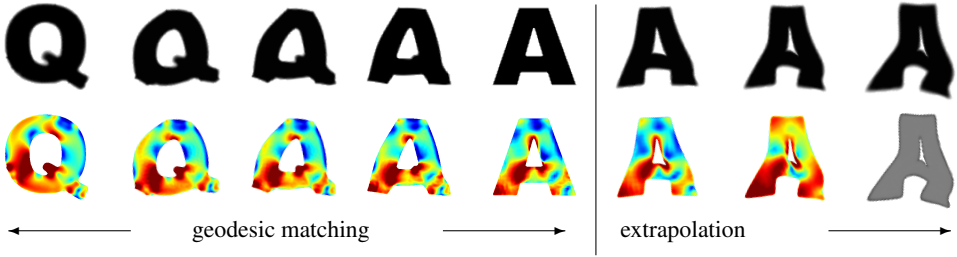

Let us introduce the space $\manifold$ of shapes as the set of all objects $\object$
which are closed subsets of $\R^d$ ($d=2,3$) and homeomorphic to a given regular reference object
$\refobject$, \ie $\object = \phi(\refobject)$ for an orientation preserving homeomorphism $\phi$.
Furthermore, objects which coincide up to a rigid body motion are identified with each other.
A smooth path $(\object(t))_{t\in [0,1]}$ in this shape space is associated with a smooth family
$(\phi(t))_{t\in [0,1]}$ of deformations.
To measure the distance between two objects, a Riemannian metric is defined on variations $\delta \object$ of
objects $\object\in\manifold$ which reflects the internal fluid friction --- called dissipation --- that occurs during the shape variation.
The local temporal rate of dissipation in a fluid depends on the symmetric part $\epsilon[v]:= \frac12(\nabla v +  \nabla v^T)$
of the gradient of the fluid velocity $v:\object\to\R^d$ (the antisymmetric remainder reflects infinitesimal rotations),
and for an isotropic Newtonian fluid,
we obtain the local rate of dissipation
\begin{equation}\label{eqn:localDissipationRate}
\dissdensity(\nabla v)= \lambda(\tr\epsilon[v])^2+2\mu\tr(\epsilon[v]^2)\,,
\end{equation}
where $\lambda,\mu$ are material-specific parameters.
Given a family $(\phi(t))_{t\in [0,1]}$ of deformations of the reference object $\refobject$, 
the change of shape along the path $(\object(t))_{t\in[0,1]}$ can be
described by the (Lagrangian) temporal variation $\dot \phi(t)$ or the associated (Eulerian) velocity field
$$v(t) = \dot \phi(t) \circ \phi^{-1}(t)$$ on $\object$. Hence, the tangent space $T_\object\manifold$ to $\manifold$ at a shape $\object$ can be identified with the space of initial velocities $v = \dot \phi(0) \circ \phi^{-1}(0)$  for deformation paths with $\phi(0,\refobject)=\object$.
Here we identify those
velocities $v$ which lead to the same effective shape variation, \ie those with the same normal component $v \cdot n$ on $\partial \object$,  where $n$ is the outer normal on $\partial \object$.
Now, integrating the local rate of dissipation for velocity fields $v$ on $\object=\phi(0,\refobject)$, we define the Riemannian metric $\metric_{\object}$ on $T_\object\manifold$ as the symmetric quadratic form with
\begin{equation}\label{eqn:viscousMetric}
\metric_{\object}(v,v)=\min\limits_{\setof{\tilde v}{
\tilde v\cdot n= v \cdot n\,\mathrm{on}\,\partial\object}} \int_\object \dissdensity(\nabla \tilde v(x))\,\d x\,.
\end{equation}
For the shape variation along a path $\curve:[0,1]\to\manifold$ described by the Eulerian motion field $(v(t))_{t\in[0,1]}$, path length $\pathlength$ and energy $\pathenergy$ are defined as
\begin{align}
\pathlength[(\curve(t))_{t\in[0,1]}]&=\textstyle\int_0^1\sqrt{\metric_{\curve(t)}(v(t),v(t))}\,\d t\,,\label{eqn:pathLength}\\
\pathenergy[(\curve(t))_{t\in[0,1]}]&=\textstyle\int_0^1\metric_{\curve(t)}(v(t),v(t))\,\d t\,.\label{eqn:pathEnergy}
\end{align}
Paths which (locally) minimize the energy $\pathenergy$ or equivalently the length $\pathlength$ are called geodesics
(\conf Figure~\ref{fig:QA}).
A geodesic thus mimics the energetically optimal way to continuously deform a fluid volume.
\subsection{Approximating the distance}
\label{sec:approxdistance}
The evaluation of the geodesic distance based on a direct space and time discretization of \eqref{eqn:viscousMetric} and \eqref{eqn:pathLength} turns out to be computationally very demanding (\conf for instance the approaches in \cite{BeMiTr05,DuGrMi98}). Hence, we use here an efficient and robust time discrete approximation based on an energy functional $\W$ which locally behaves like the squared Riemannian distance (\ie the squared length of a connecting geodesic): \medskip

Given two shapes $\object$ and $\tilde \object$, we consider an approximation
\begin{equation}
\label{eq:approximation}
\dist^2(\object,\tilde \object) \approx \W_\object[\psi]\,,
\end{equation}
where $\W_\object[\tilde\psi] := \int_\object W(\nabla\tilde\psi) \,\d x$ is the stored deformation energy of a deformation $\tilde\psi:\object\to\R^d$
and $\psi$ is the minimizer of this energy over all such deformations with $\psi(\object)=\tilde\object$.
Here, $W:\R^{d,d}\to[0,\infty)$ is a so-called hyperelastic energy density.\medskip

In correspondence to our assumption that objects are identical if they coincide up to a rigid body motion,
we require $W$ to be rigid body motion invariant. Furthermore, we assume the objects to have no preferred material directions so that
$W$ is in addition isotropic, which altogether leads to $W(RAU)=W(A)$  for all $R,U\in SO(d),A\in\R^{d,d}$ (\conf \cite{Ci97}).
In the undeformed configuration for $\psi=\id$, energy and stresses (the first derivatives of $W$) are supposed to vanish so that we require $W(\Id)=0$, $\D W(\Id)=0$ (where $\D W$ denotes the derivative with respect to the matrix argument).
Furthermore, we need $W(A)\to\infty$ as $\det A\to0$ to prohibit material self-penetration, which is linked to the preservation of topology.
The approximation property \eqref{eq:approximation} relies on a consistent choice of $\W_\object$ for the given metric $\metric_{\object}$
which can be expressed by the relation
\begin{equation}\label{eqn:energyHessianIsMetric}
\frac12\left.\frac{\mathrm{d}^2}{\mathrm{d}t^2} \W_{\object}[\psi(t)]\right|_{t=0} = \frac12 \int_{\object} \D^{2} W(\Id)(\nabla v,\nabla v)\, \d x
= \int_{\object} \dissdensity(\nabla v) \,\d x
\end{equation}
along any object path $\object(t)=\psi(t,\object)$, $t\in\R$, with $\psi(0)= \id$ and  velocity field $v=\dot \psi(0)$.
Using the notion of the Hessian $\mathrm{Hess}_\manifold$ of  a function on a manifold as the
endomorphism representing its second variation in the metric, we can rephrase this approximation condition more geometrically as $$\frac12 \mathrm{Hess}_\manifold \W_\object[\id] = \id$$ with the usual identification of objects $\object$ and deformations $\phi$.
For the deformation energy density $W$, this condition implies that
its Hessian $\D^2 W(\Id)$ has to satisfy $ \frac12 \D^2 W(\Id)(A,A) = \dissdensity(A)$ for all $A \in \R^{d,d}\,$.
A suitable example is
\begin{equation*}
W(A)=\frac\mu2\tr (A^TA)+\frac\lambda4\det A^2-\left(\mu+\frac\lambda2\right)\log\det A-\frac{d\mu}2-\frac\lambda4\,.
\end{equation*}
Assume that the energy density satisfies the above-mentioned properties. We observe that
the metric $\metric_{\object}$ is the first non-vanishing term in the Taylor expansion of the squared length of a curve, \ie
$$\left(\pathlength[(\curve(t))_{t\in[0,T]}]\right)^2 = T^2\metric_{\curve(0)}(v,v) + O(T^3)$$
with $v=\dot \phi(0)\circ \phi^{-1}(0)$ being the initial tangent vector along a smooth path $(\curve(t))_{t\in[0,T]}=(\phi(t,\refobject))_{t\in [0,1]}$. Thus, since the Hessian of the energy $\W_\object$ and the metric $\metric_\object$ are related by \eqref{eqn:energyHessianIsMetric}, we obtain that the second order Taylor expansions of $\dist^2(\object,\psi(\object))$ and $\W_\object[\psi]$ in $\psi$ coincide and indeed
\begin{equation}
\label{eq:distapprox}
\dist^2(\object,\tilde \object) = \min_{\setof{\psi}{\psi(\object)=\tilde \object}} \W_\object[\psi] + O(\dist^3(\object,\tilde \object))\,.
\end{equation}
Here, different from \cite{WiBaRu10} we neither take into account mismatch penalties nor perimeter regularizing functionals for each object $\object_k$, $k=0,\ldots, K$.

\subsection{Discrete length and discrete energy}
\label{sec:discretelengthenergy}
Now, we are in a position to discretize length and energy of paths $(\object(t))_{t\in [0,1]}$ in shape space.
To this end, we first sample the path at times $t_k=k \tau$ for $k=0 ,\ldots, K$ ($\tau= \frac1K$), denote $\object_k:=\object(t_k)$,
and obtain the estimates 
\begin{eqnarray*}
\pathlength[(\curve(t))_{t\in[0,1]}] &\geq&\textstyle \sum_{k=1}^K\dist(\object_{k-1},\object_k) \\
\pathenergy[(\curve(t))_{t\in[0,1]}]&\geq&\textstyle \frac1\tau\sum_{k=1}^K\dist^2(\object_{k-1},\object_k)
\end{eqnarray*}
for the length and the energy, where equality holds for geodesic paths. Indeed, the first estimate is straightforward, and the application of the Cauchy--Schwarz inequality leads to
\begin{eqnarray*}
\sum_{k=1}^K\dist^2(\object_{k-1},\object_k) 
&\leq& \sum_{k=1}^K \left(\int_{(k-1)\tau}^{k\tau} \sqrt{\metric_{\curve(t)}(v(t),v(t))}\, \d t\right)^2 \\
&\leq& \sum_{k=1}^K \tau \; \int_{(k-1)\tau}^{k\tau} \metric_{\curve(t)}(v(t),v(t))\, \d t = \tau \; \pathenergy[(\curve(t))_{t\in[0,1]}]\\
\end{eqnarray*}
which implies the second estimate.

Together with \eqref{eq:distapprox} this motivates the following definition of a \emph{discrete path energy} and a \emph{discrete path length}
for a discrete path $(\object_0,\ldots,\object_K)$ in shape space:
\begin{align}
\pathlength[(\object_0,\ldots,\object_K)]&=\textstyle\sum_{k=1}^K \sqrt{\W_{\curve_{k-1}}[\psi_k]}\,,\label{eqn:discretePathLength}\\
\pathenergy[(\object_0,\ldots,\object_K)]&=\textstyle \frac1\tau\sum_{k=1}^K \W_{\object_{k-1}}[\psi_k]\,, \label{eqn:discretePathEnergy}
\end{align}
where $\psi_k = \argmin_{\setof{\psi}{\psi(\object_{k-1})= \object_k}} \W_{\object_{k-1}}[\psi]$ (\conf also \cite{WiBaRu10}).
In fact, \eqref{eqn:discretePathLength} and \eqref{eqn:discretePathEnergy} can for general smooth paths even be proven to be first order consistent with
the continuous length \eqref{eqn:pathLength} and energy \eqref{eqn:pathEnergy} as $\tau\to0$.
For illustration, if $\manifold$ is a two-dimensional manifold embedded in $\R^3$, we can interpret the terms $\W_{\object_{k-1}}$ as the stored elastic energies in springs which connect a sequence of points $\object_k$ on the manifold through the ambient space.
Then the discrete path energy is the total stored elastic energy in this chain of springs.

\begin{figure*}[t]
\centering
\setlength{\unitlength}{.112\linewidth}
\begin{picture}(8.86,1.86)
\put(0,1){\includegraphics[width=.86\unitlength]{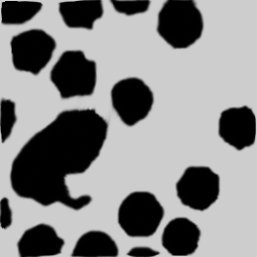}}
\put(1,1){\includegraphics[width=.86\unitlength]{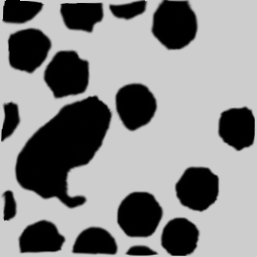}}
\put(2,1){\includegraphics[width=.86\unitlength]{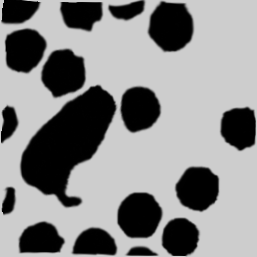}}
\put(3,1){\includegraphics[width=.86\unitlength]{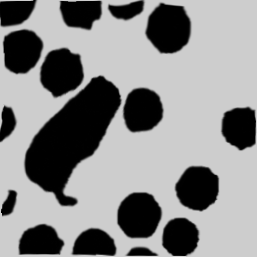}}
\put(4,1){\includegraphics[width=.86\unitlength]{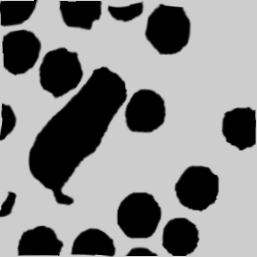}}
\put(5,1){\includegraphics[width=.86\unitlength]{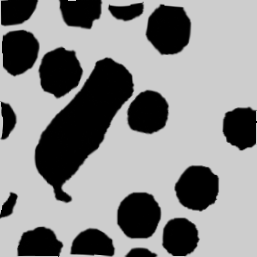}}
\put(6,1){\includegraphics[width=.86\unitlength]{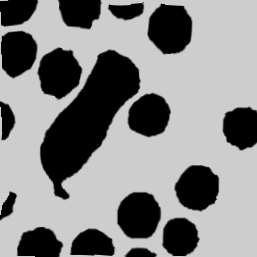}}
\put(7,1){\includegraphics[width=.86\unitlength]{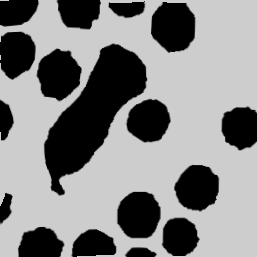}}
\put(8,1){\includegraphics[width=.86\unitlength]{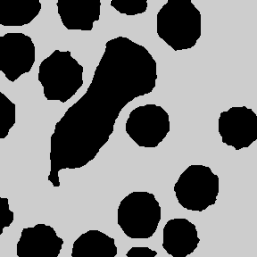}}
\put(0,0){\includegraphics[width=.86\unitlength]{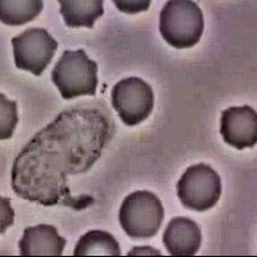}}
\put(1,0){\includegraphics[width=.86\unitlength]{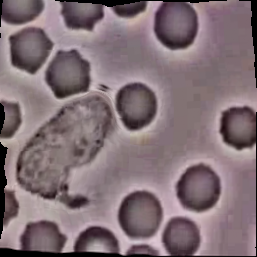}}
\put(2,0){\includegraphics[width=.86\unitlength]{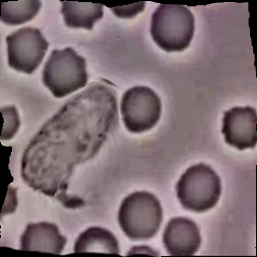}}
\put(3,0){\includegraphics[width=.86\unitlength]{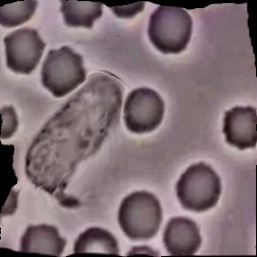}}
\put(4,0){\includegraphics[width=.86\unitlength]{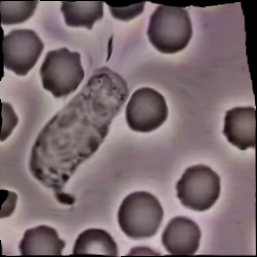}}
\put(5,0){\includegraphics[width=.86\unitlength]{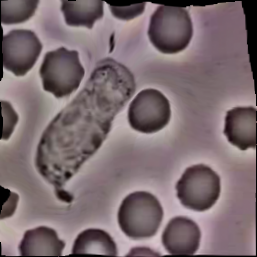}}
\put(6,0){\includegraphics[width=.86\unitlength]{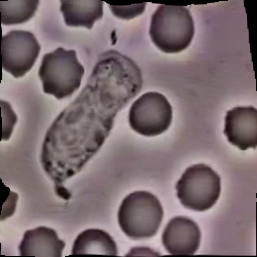}}
\put(7,0){\includegraphics[width=.86\unitlength]{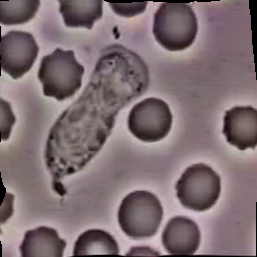}}
\put(8,0){\includegraphics[width=.86\unitlength]{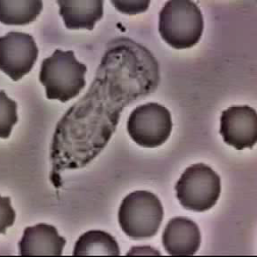}}
\put(.88,.03){\line(1,4){.1}}
\put(.88,.83){\line(1,-4){.1}}
\put(1.88,.03){\line(1,4){.1}}
\put(1.88,.83){\line(1,-4){.1}}
\put(2.88,.03){\line(1,4){.1}}
\put(2.88,.83){\line(1,-4){.1}}
\put(5,.03){\line(-1,4){.1}}
\put(5,.83){\line(-1,-4){.1}}
\put(6,.03){\line(-1,4){.1}}
\put(6,.83){\line(-1,-4){.1}}
\put(7,.03){\line(-1,4){.1}}
\put(7,.83){\line(-1,-4){.1}}
\put(8,.03){\line(-1,4){.1}}
\put(8,.83){\line(-1,-4){.1}}
\end{picture}
\caption{Nonlinear video interpolation via a discrete geodesic (top) between two segmented photographs of white and red blood cells
(first and last picture of bottom row, courtesy Robert A.\ Freitas, Institute for Molecular Manufacturing, California, USA).
The bottom row shows pushforwards and pullbacks of the end images under the deformations along the discrete geodesic.}
\label{fig:bloodCells}
\end{figure*}
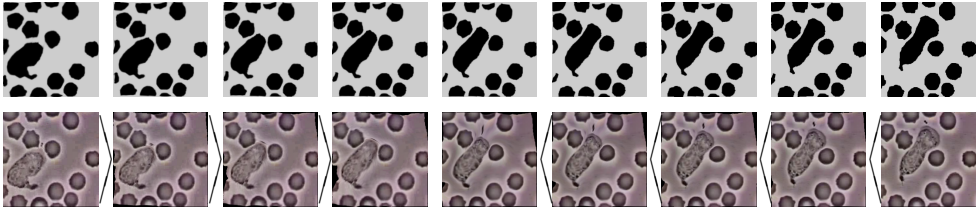
\medskip

A \emph{discrete geodesic} (of order $K$) is now defined
as a minimizer of $\pathenergy[(\curve_0,\ldots,\curve_K)]$ for fixed end points $\curve_0,\curve_K$. 
The discrete geodesic is thus an energetically optimal sequence of deformations from $\object_0$ into $\object_K$.
\medskip

In the minimization algorithm to be discussed in Section\,\ref{sec:parageodesic} we do not explicitly minimize $\pathenergy[(\curve_0,\ldots,\curve_K)]$ for the object contours as in \cite{WiBaRu10} but instead for reference deformations defined on fixed reference objects.  Figure~\ref{fig:bloodCells} shows a discrete geodesic in the context of multicomponent objects, which is visually identical to that obtained by the more complex approach in \cite{WiBaRu10}. Here, deformations are considered which map every component of a shape onto the corresponding component of the next shape in the discrete path as the obvious generalization of discrete geodesics between single component shapes.


\begin{wrapfigure}{r}{0.3\linewidth}
\centering\includegraphics[width=0.95\linewidth,height=.4\linewidth]{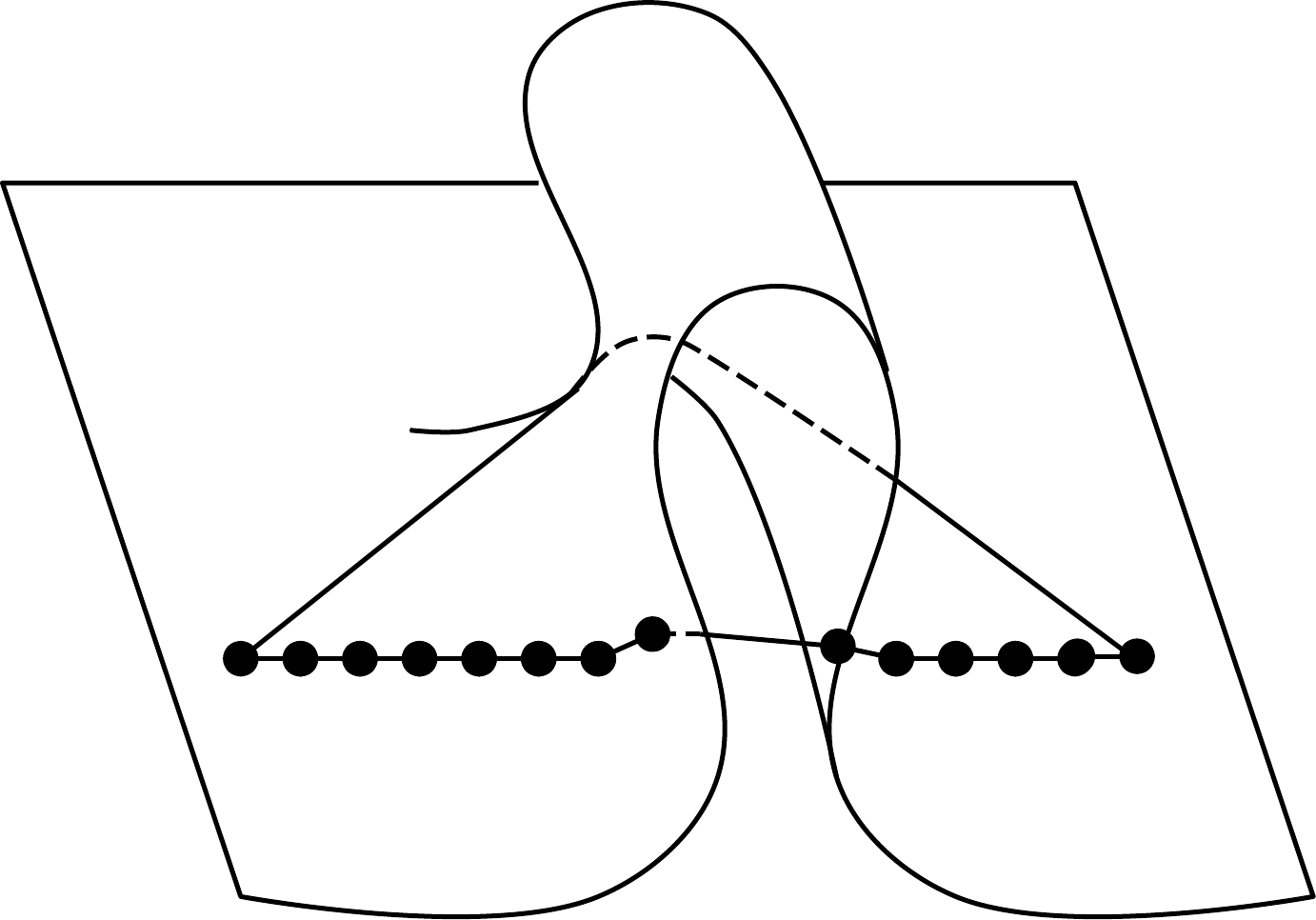}
\caption{A continuous geodesic and a discrete path which almost minimizes the discrete path length on a two-dimensional manifold embedded in $\R^3$.
}
\label{fig:energyvslength}
\end{wrapfigure}
While in the continuous case geodesic curves equally minimize length \eqref{eqn:pathLength} and energy \eqref{eqn:pathEnergy},
minimizers of the discrete path length \eqref{eqn:discretePathEnergy} are in general not related to discrete geodesics (and thus also not to continuous geodesics as $\tau\to0$).
Indeed, let us consider a two-dimensional manifold $\manifold$ embedded in $\R^3$,
paired with the deformation energy $\W_{\object_{k-1}}[\id + \zeta_k]:=|\zeta_k|^2$ for a displacement vector 
$\zeta_k$ in $\R^3$
connecting points $\object_{k-1}$ and $\object_{k}$ on $\manifold$.
Now take into account a continuous geodesic and a discrete path on $\manifold$ where the end points are close to each other in the embedding space but far apart on the surface. Figure~\ref{fig:energyvslength} depicts such a configuration with a discrete path which almost minimizes the discrete path length.
A minimizer of the discrete path length will always jump through the protrusion and never approximate the continuous geodesic, whereas minimizers of the discrete path energy satisfy $\W_{\object_{k-1}}[\id +\zeta_k]\to 0$ as $\tau \to 0$ and thus rule out such a short cut through the ambient space.

\section{Time discrete geodesic calculus}
\label{sec:timediscrete}
With the notion of discrete geodesics at hand we will now derive a full-fledged discrete geodesic calculus based on a time discrete geometric logarithm and a time discrete exponential map, which then also give rise to a discrete parallel transport and a discrete
Levi-Civita connection on shape space.

\subsection{Discrete logarithm and shape variations}
If $(\object(t))_{t\in[0,1]}$ is the unique geodesic on $\manifold$ connecting $\object=\object(0)$ and $\tilde \object=\object(1)$,
the logarithm of $\tilde \object$ with respect to $\object$ is defined as the initial velocity $v \in T_\object\manifold$ of the geodesic
path. In terms of Section\,\ref{sec:viscousshapespace} we have
\begin{equation*}
\log_\object(\tilde \object) =v(0)
\end{equation*}
for $v(t)=\dot\phi(t)\circ \phi(t)^{-1}$, where $\phi(t,\refobject) = \object(t)$ defines the associated family of deformations.
On a geodesically complete Riemannian manifold the logarithm exists as long as $\dist(\object, \tilde \object)$ is sufficiently small. The associated logarithmic map
$\log_\object: \tilde \object \mapsto v(0)\in T_\object\manifold$ represents (nonlinear) variations on the manifold as (linear) tangent vectors.

\begin{figure}
\setlength{\unitlength}{.14\linewidth}
\hspace{0.3\linewidth}
\begin{picture}(5,3.9)(.1,-1.4)
\put(-1.4,2.55){computing}
\put(-1.4,2.3){$K\Log{K}_{P}(A)$:}
\put(2.7,.1){\includegraphics[width=1.2\unitlength]{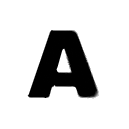}}
\put(4.1,.7){\includegraphics[width=.4\unitlength]{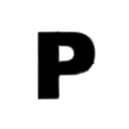}}
\put(4.8,.7){\includegraphics[width=.4\unitlength]{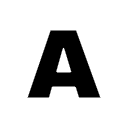}}
\put(2.0,.5){\includegraphics[width=1.2\unitlength,angle=0,origin=c]{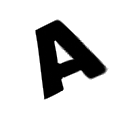}}
\put(3.4,1.1){\includegraphics[width=.4\unitlength]{PA1}}
\put(4.1,1.1){\includegraphics[width=.4\unitlength]{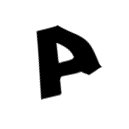}}
\put(4.8,1.1){\includegraphics[width=.4\unitlength]{PA5}}
\put(1.2,0.9){\includegraphics[width=1.2\unitlength,angle=0,origin=c]{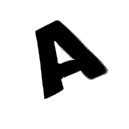}}
\put(2.7,1.5){\includegraphics[width=.4\unitlength]{PA1}}
\put(3.4,1.5){\includegraphics[width=.4\unitlength]{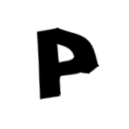}}
\put(4.1,1.5){\includegraphics[width=.4\unitlength]{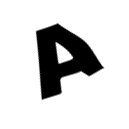}}
\put(4.8,1.5){\includegraphics[width=.4\unitlength]{PA5}}
\put(.4,1.3){\includegraphics[width=1.2\unitlength,angle=0,origin=c]{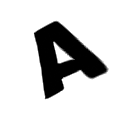}}
\put(2.0,1.9){\includegraphics[width=.4\unitlength]{PA1}}
\put(2.7,1.9){\includegraphics[width=.4\unitlength]{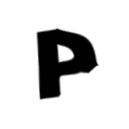}}
\put(3.4,1.9){\includegraphics[width=.4\unitlength]{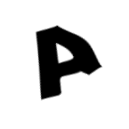}}
\put(4.1,1.9){\includegraphics[width=.4\unitlength]{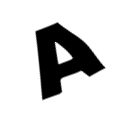}}
\put(4.8,1.9){\includegraphics[width=.4\unitlength]{PA5}}
\put(-.3,1.7){\includegraphics[width=1.2\unitlength]{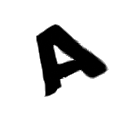}}
\put(1.3,2.3){\includegraphics[width=.4\unitlength]{PA1}}
\put(2.0,2.3){\includegraphics[width=.4\unitlength]{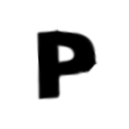}}
\put(2.7,2.3){\includegraphics[width=.4\unitlength]{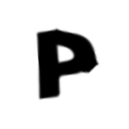}}
\put(3.4,2.3){\includegraphics[width=.4\unitlength]{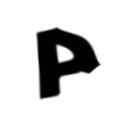}}
\put(4.8,2.3){\includegraphics[width=.4\unitlength]{PA5}}
\put(4.45,0.85){\vector(1,0){.45}}
\put(4.15,0.85){\vector(-1,0){.6}}
\put(3.425,0.975){\tiny$1(\!\frac11\!\mathrm{LOG}\!)_{\!P}\!(\!A\!)$}
\put(4.45,1.25){\vector(1,0){.45}}
\put(3.75,1.25){\vector(1,0){.45}}
\put(3.45,1.25){\vector(-1,0){.6}}
\put(2.725,1.375){\tiny$2(\!\frac12\!\mathrm{LOG}\!)_{\!P}\!(\!A\!)$}
\put(4.45,1.65){\vector(1,0){.45}}
\put(3.75,1.65){\vector(1,0){.45}}
\put(3.05,1.65){\vector(1,0){.45}}
\put(2.75,1.65){\vector(-1,0){.6}}
\put(2.0,1.75){\tiny$3(\!\frac13\!\mathrm{LOG}\!)_{\!P}\!(\!A\!)$}
\put(4.45,2.05){\vector(1,0){.45}}
\put(3.75,2.05){\vector(1,0){.45}}
\put(3.05,2.05){\vector(1,0){.45}}
\put(2.35,2.05){\vector(1,0){.45}}
\put(2.05,2.05){\vector(-1,0){.6}}
\put(1.3,2.15){\tiny$4(\!\frac14\!\mathrm{LOG}\!)_{\!P}\!(\!A\!)$}
\put(4.45,2.45){\vector(1,0){.45}}
\put(3.75,2.45){\vector(1,0){.45}}
\put(3.05,2.45){\vector(1,0){.45}}
\put(2.35,2.45){\vector(1,0){.45}}
\put(1.65,2.45){\vector(1,0){.45}}
\put(1.45,2.45){\vector(-1,0){.6}}
\put(0.6,2.55){\tiny$8(\!\frac18\!\mathrm{LOG}\!)_{\!P}\!(\!A\!)$}
\put(4.2,2.45){$\ldots$}
\end{picture}
\vspace{-8.6cm}

\begin{picture}(5,3.9)(.1,-1.4)
\put(0.1,1.25){applying}
\put(0.1,1.0){$\ExpO{K}{P}\!(\frac\zeta{K})$:}
\put(0,0.2){\includegraphics[width=.4\unitlength]{PA1}}
\put(0.7,0){\includegraphics[width=1.2\unitlength]{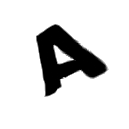}}
\put(0,-.2){\includegraphics[width=.4\unitlength]{PA1}}
\put(0.9,-.2){\includegraphics[width=.4\unitlength]{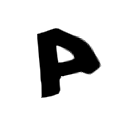}}
\put(1.6,-.5){\includegraphics[width=1.2\unitlength,angle=0,origin=c]{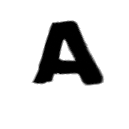}}
\put(0,-.6){\includegraphics[width=.4\unitlength]{PA1}}
\put(0.9,-.6){\includegraphics[width=.4\unitlength]{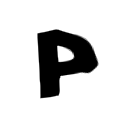}}
\put(1.8,-.6){\includegraphics[width=.4\unitlength]{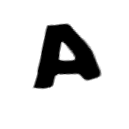}}
\put(2.5,-.9){\includegraphics[width=1.2\unitlength,angle=0,origin=c]{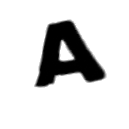}}
\put(0,-1.0){\includegraphics[width=.4\unitlength]{PA1}}
\put(0.9,-1.0){\includegraphics[width=.4\unitlength]{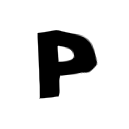}}
\put(1.8,-1.0){\includegraphics[width=.4\unitlength]{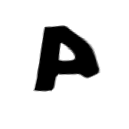}}
\put(2.7,-1.0){\includegraphics[width=.4\unitlength]{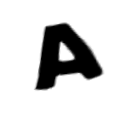}}
\put(3.4,-1.3){\includegraphics[width=1.2\unitlength,angle=0,origin=c]{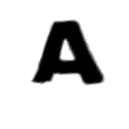}}
\put(0,-1.4){\includegraphics[width=.4\unitlength]{PA1}}
\put(0.9,-1.4){\includegraphics[width=.4\unitlength]{ExpPA8b}}
\put(1.8,-1.4){\includegraphics[width=.4\unitlength]{ExpPA8c}}
\put(2.7,-1.4){\includegraphics[width=.4\unitlength]{ExpPA8d}}
\put(4.3,-1.6){\includegraphics[width=1.2\unitlength,angle=0,origin=c]{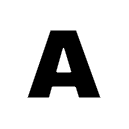}}
\put(.4,.4){\vector(1,0){.5}}
\put(.4,.47){\tiny$\id + \zeta$}
\put(.4,.0){\vector(1,0){.5}}
\put(1.3,.0){\vector(1,0){.5}}
\put(.4,.07){\tiny$\id + \zeta/2$}
\put(.4,-.4){\vector(1,0){.5}}
\put(1.3,-.4){\vector(1,0){.5}}
\put(2.2,-.4){\vector(1,0){.5}}
\put(.4,-.33){\tiny$\id + \zeta/3$}
\put(.4,-.8){\vector(1,0){.5}}
\put(1.3,-.8){\vector(1,0){.5}}
\put(2.2,-.8){\vector(1,0){.5}}
\put(3.1,-.8){\vector(1,0){.5}}
\put(.4,-.73){\tiny$\id + \zeta/4$}
\put(.4,-1.2){\vector(1,0){.5}}
\put(1.3,-1.2){\vector(1,0){.5}}
\put(2.2,-1.2){\vector(1,0){.5}}
\put(3.1,-1.2){\vector(1,0){.5}}
\put(4.0,-1.2){\vector(1,0){.5}}
\put(.4,-1.13){\tiny$\id + \zeta/8$}
\put(3.7,-1.2){$\ldots$}
\end{picture}
\caption{In the bottom left we experimentally verify the convergence stated in \eqref{eq:EXPconvergence}
by computing discrete geodesics starting from `P' and ending at $\ExpO{K}{P}(\zeta/K)$ for $K=1,2,3,4,8$ and $\zeta=8\Log{8}_P(A)$.
In the top right, based on a computation of discrete geodesics between `P' and `A' of order $K=8, 4, 3, 2, 1$, the resulting logarithms $\Log{K}_P(A)$ are depicted, where $\Log{K}_P(A)$ is the displacement from `P' to the second shape of the respective geodesic. To the left of the initial shape `P', those shapes are displayed which result from applying the displacement $K\Log{K}_P(A)$ to `P' to experimentally verify convergence as stated in \eqref{eq:LOGconvergence}.
The arbitrary rotations are due to the rigid body motion invariance of our approach.
}
\label{fig:LogExp}
\end{figure}

The initial velocity $v(0)$ can be approximated by a difference quotient in time,
$$v(0,x)= \frac1\tau \zeta(x) + O(\tau)\,,$$
where $\zeta(x) = \phi(\tau,x)\circ \phi(0,x)^{-1}- x$ denotes a displacement on the initial object $\object$.
Thus, we obtain  $$\tau \log_\object(\tilde \object) = \zeta(x)+ O(\tau^2)\,.$$
This gives rise to a consistent definition of a time discrete logarithm. Let $(\object_0,\ldots,\object_K)$ be a discrete geodesic
between $\object=\object_0$ and $\tilde\object=\object_K$
with an associated sequence of matching deformations $\psi_1,\ldots, \psi_K$, then
we consider $\frac1\tau \zeta_1$ for the displacement $\zeta_1(x) = \psi_1(x)- x$
as an approximation of $v(0)=\log_\object(\tilde \object)$. Taking into account that $\tau = \frac1K$ we thus define the discrete
$\frac1K$-logarithm
\begin{equation}\label{eq:LOGK}
\textstyle\Log{K}_\object(\tilde\object) := \zeta_1\,.
\end{equation}
In the special case $K=1$ and a discrete geodesic $(\object,\tilde \object)$
we simply obtain
\begin{equation*}
\Log{1}_\object(\tilde\object) =
\argmin_{\setof{\zeta_1}{(\id+\zeta_1)(\object) = \tilde \object}} \W_\object[\id+\zeta_1]\,.
\end{equation*}
As in the continuous case the discrete logarithm can be considered as a
representation of the nonlinear variation $\tilde \object$ of $\object$ in the (linear) tangent space of displacements on $\object$.
On a sequence of successively refined discrete geodesics we expect
\begin{equation}\label{eq:LOGconvergence}
K \textstyle\Log{K}_\object(\tilde\object) \to \log_\object(\tilde\object)
\end{equation}
for $K\to \infty$ (\conf Figure~\ref{fig:LogExp} for an experimental validation of this convergence behaviour).

\subsection{Discrete exponential and shape extrapolation}
\label{sec:exp}
In the continuous setting, the exponential map $\exp_\object$ maps tangent vectors $v \in T_\object\manifold$  onto
the end point $\object(1)$ of a geodesic $(\object(t))_{t\in[0,1]}$ with $\object(0)=\object$ and the given tangent vector $v$ at time $0$.
Using the notation from the previous section we have $\exp_\object(v(0))= \tilde \object$ and, via a simple scaling argument,
$\textstyle\exp_\object\left(\frac{k}{K} v(0)\right) = \object(\frac{k}{K})$
for $k=0,\ldots, K$.
We now aim at defining a discrete power $k$ exponential map $\ExpO{k}{\object}$ such that $\ExpO{k}{\object}(\zeta_1)=\object_k$
on a discrete geodesic $(\object,\object_1,\ldots,\object_K)$ of order $K\geq k$ with $\zeta_1 = \Log{K}_\object(\tilde \object)$
(the notation is motivated by the observation that $\exp(ks)= \exp^k(s)$ on $\R$ or more general matrix groups).
Our definition will reflect the following recursive properties of the continuous exponential map,
\begin{align*}
\exp_\object(1v) =& \left(\textstyle\frac11\log_\object\right)^{-1}(v)\,,\\
\exp_\object(2v) =& \left(\textstyle\frac12 \log_\object\right)^{-1}(v)
\notinclude{\Leftrightarrow\textstyle\frac12 \log_\object (\exp_\object (2v)) = v}\,, \\
\exp_\object(kv) =& \exp_{\exp_\object((k-2)v)}(2 v_{k-1}) \\
&\text{for }v_{k-1} := \log_{\exp_\object ((k-2) v)} \exp_\object((k-1) v)\,.
\end{align*}
Replacing $\exp(k\cdot)$ by $\Exp{k}$, $\frac11 \log$ by $\Log{1}$, and $\frac12 \log$ by $\Log{2}$ we obtain the recursive definition
\begin{align}
\ExpO{1}{\object}(\zeta) :=& \textstyle\Log{1}_\object^{-1}(\zeta)\,, \label{eq:exp1}\\
\ExpO{2}{\object}(\zeta) :=& \textstyle\Log{2}_\object^{-1}(\zeta)\,, \label{eq:exp2}\\
\ExpO{k}{\object}(\zeta) :=& \ExpO{2}{\ExpO{k-2}{\object}(\zeta)}(\zeta_{k-1}) \label{eq:expk}\\
& \text{with }\zeta_{k-1} := \textstyle\Log{1}_{\ExpO{k-2}{\object}(\zeta)} \ExpO{k-1}{\object}(\zeta)\,. \nonumber
\end{align}
It is straightforward to verify that $\ExpO{K}{\object}=\Log{K}_\object^{-1}$ as long as the discrete logarithm on the right is invertible.
Equation\,\eqref{eq:exp1} implies $\ExpO{1}{\object}(\zeta)=(\id+\zeta)(\object)$,
and \eqref{eq:exp2} in fact represents a variational constraint for a discrete geodesic flow of order $2$\,:
\medskip

{\it Given the object $\object$ we consider discrete geodesic paths $(\object, \tilde \object_1, \tilde \object_2)$ of order $2$,
where for any chosen $\tilde \object_2$ the middle object $\tilde\object_1$ is defined via minimization of \eqref{eqn:discretePathEnergy}
so that we may write $\tilde \object_1[\tilde \object_2]$.
We now identify $\ExpO{2}{\object} (\zeta)$ as the object $\tilde \object_2$ for which $(\id+\zeta)(\object)=\tilde\object_1[\tilde\object_2]$,
\ie $\zeta$ is the energetically optimal displacement from $\object$ to $\tilde \object_1[\tilde \object_2]$ and thus satisfies
\begin{equation}\label{eqn:variationalProblemExp2}
\id + \zeta = \argmin\limits_{\setof{\psi_1}{\psi_1(\object)=\tilde \object_1[\tilde \object_2]}} \W_{\object}[\psi_1]
\end{equation}
up to a rigid body motion.
}
\medskip

Alternatively, the condition \eqref{eqn:variationalProblemExp2} can be phrased as 
\begin{equation}\label{eqn:altvariationalProblemExp2}
\id + \zeta = \argmin_{\{\psi_1\}} \mathrm{min}_{\setof{\psi_2}{(\psi_2\circ\psi_1)(\object)=\tilde \object_2}} 
\left(\W_{\object}[\psi_1] + \W_{\psi_1(\object)}[\psi_2]\right)\,.
\end{equation}

\begin{wrapfigure}{r}{0.5\linewidth}
\setlength{\unitlength}{\linewidth}
\begin{picture}(1,.35)(0,.05)
\put(0,.05){\includegraphics[width=.9\unitlength,angle=-20,origin=c]{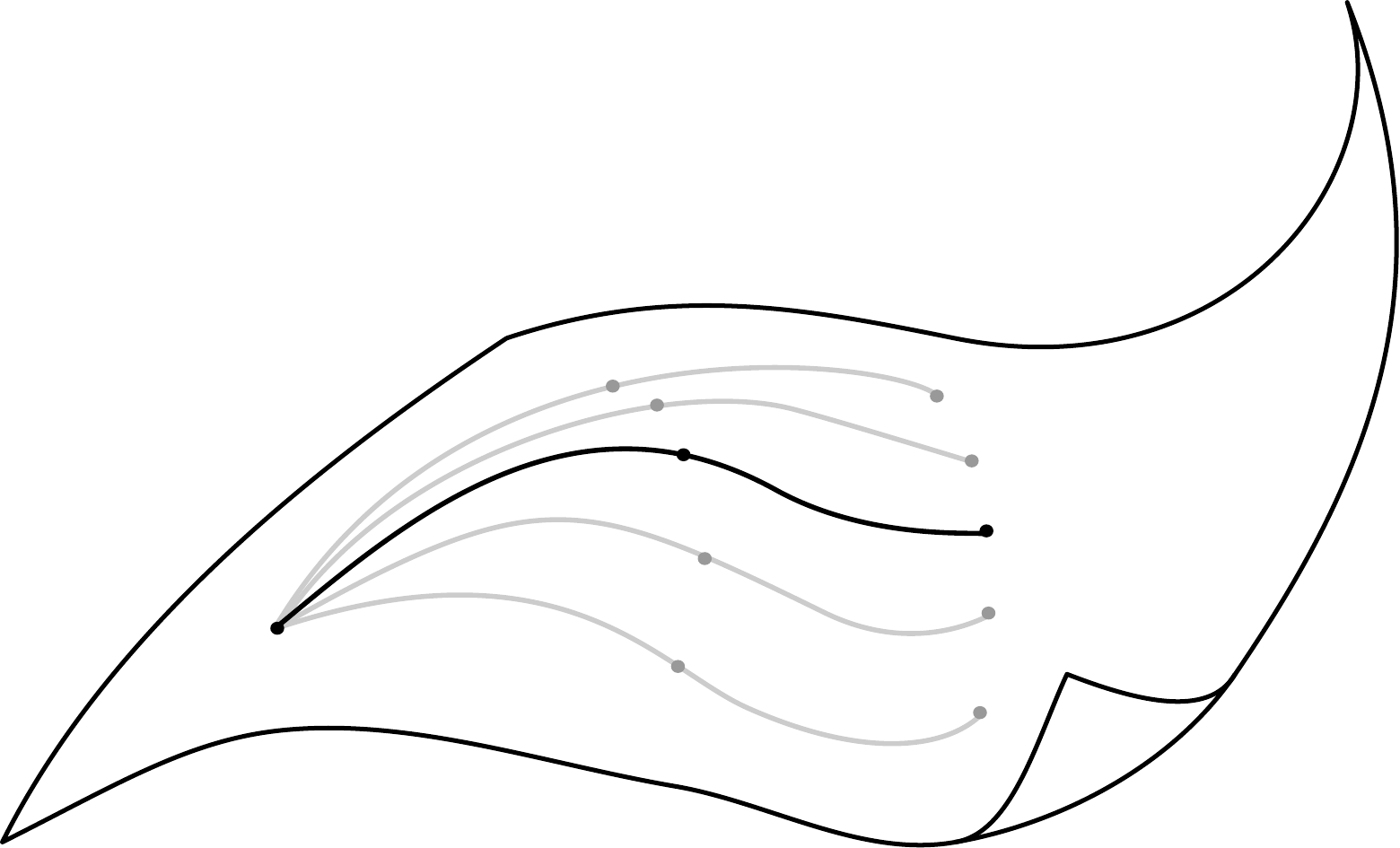}}
\put(.1,.19){\small$\manifold$}
\put(.29,.27){\small$(\id\!+\!\zeta)(\!\object\!)$}
\put(.16,.27){\small$\object$}
\put(.68,.18){\small$\ExpO{2}{\object}(\zeta)$}
\put(.62,.09){\small$\tilde{\object}_2$}
\put(.39,.135){\small$\tilde{\object}_1[\tilde{\object}_2]$}
\end{picture}
\caption{Conceptual sketch of the procedure to compute $\ExpO{2}{\object}(\zeta)$ with objects represented as points.}
\label{fig:sketchExp}
\end{wrapfigure}
Figure~\ref{fig:sketchExp} conceptually sketches the procedure to
compute $\ExpO{2}{\object}(\zeta)$. 
For given initial object $\object$ and initial displacement $\zeta$ the discrete exponential
$\ExpO{2}{\object} (\zeta)$ is selected from a fan of discrete geodesics with varying $\tilde \object_2$ as the terminal point of a discrete geodesic of order $2$ in such a way that \eqref{eqn:altvariationalProblemExp2} holds.
To compute $\ExpO{2}{\object} (\zeta)$ in the geodesic flow
algorithm \eqref{eq:exp2} and \eqref{eq:expk} we have to find the root of
\begin{equation}\label{eqn:fixedPointExp2}
F_{\object,\zeta}(\tilde \object_2) =  \Log{2}_\object(\tilde \object_2)-\zeta\,,
\end{equation}
implicitly assuming that $\zeta$ is small enough so that discrete geodesics are unique (\conf Section\,\ref{sec:Exp2numerics} for the algorithmic realization based on a representation of the unknown domain $\tilde \object_2$ via a deformation).

Equation\,\eqref{eq:expk} describes the recursion to compute $\ExpO{k}{\object}$ based on the above $\ExpO{2}{}$ single step scheme: For given $\object_{k-2}=\ExpO{k-2}{\object}(\zeta)$ and
$\object_{k-1}=\ExpO{k-1}{\object}(\zeta)$ one first retrieves $\zeta_{k-1}=\psi_{k-1}-\id$ from the previous step, where
$$\psi_{k-1}=\argmin_{\setof{\psi}{\psi(\object_{k-2})=\object_{k-1}}} \W_{\object_{k-2}}[\psi]\,.$$
Then \eqref{eq:exp2} is applied to compute $\object_k$ from $\object_{k-2}$ and $\zeta_{k-1}$ as the root of
$F_{\object_{k-2},\zeta_{k-1}}$.

For sufficiently small $\zeta$ we expect
$\Exp{k}$ to be well-defined.
Since by definition, every triplet $(\object_{k-1},\object_k,\object_{k+1})$ of the sequence $\object_k=\ExpO{k}{\object}(\zeta)$ is a geodesic of order 2 and minimizes $\pathenergy[(\object_{k-1},\object_k,\object_{k+1})]$, the resulting family $(\object_k)_{k=0,\ldots, K}$ indeed is a discrete geodesic of order $K$.
\begin{figure}
\begin{center}
\setlength{\unitlength}{1.2\linewidth}
\begin{picture}(.63,.2)
\put(-.01,.02){\includegraphics[width=.2\unitlength]{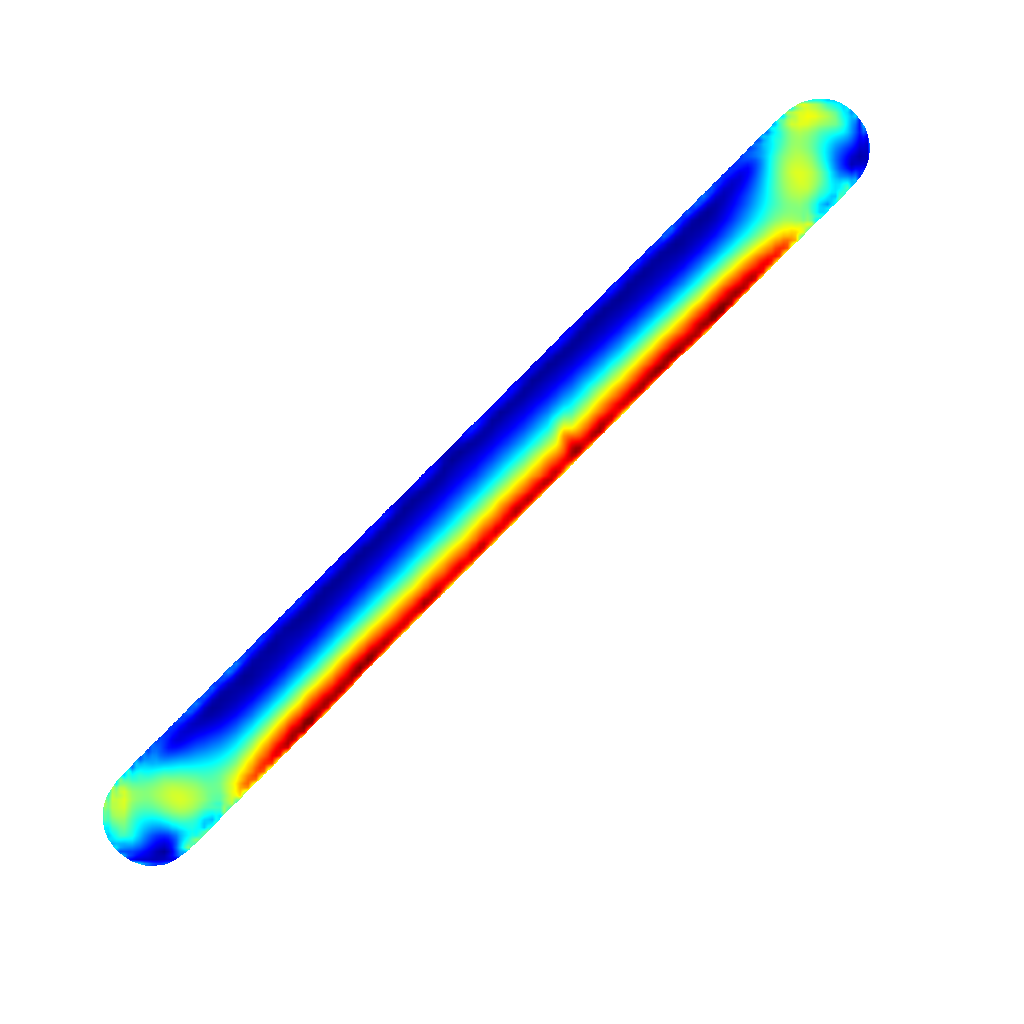}}
\put(0.03,.01){\includegraphics[width=.2\unitlength]{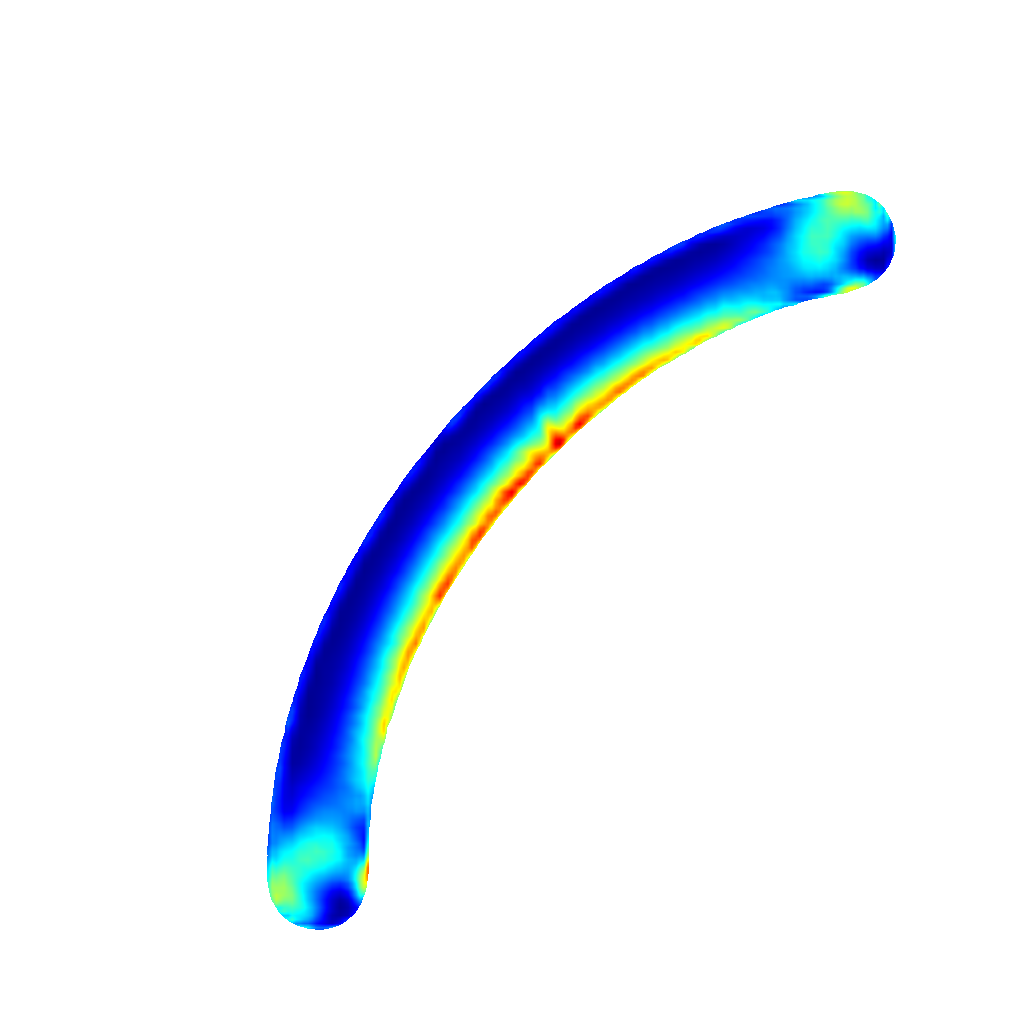}}
\put(.07,-.01){\includegraphics[width=.2\unitlength]{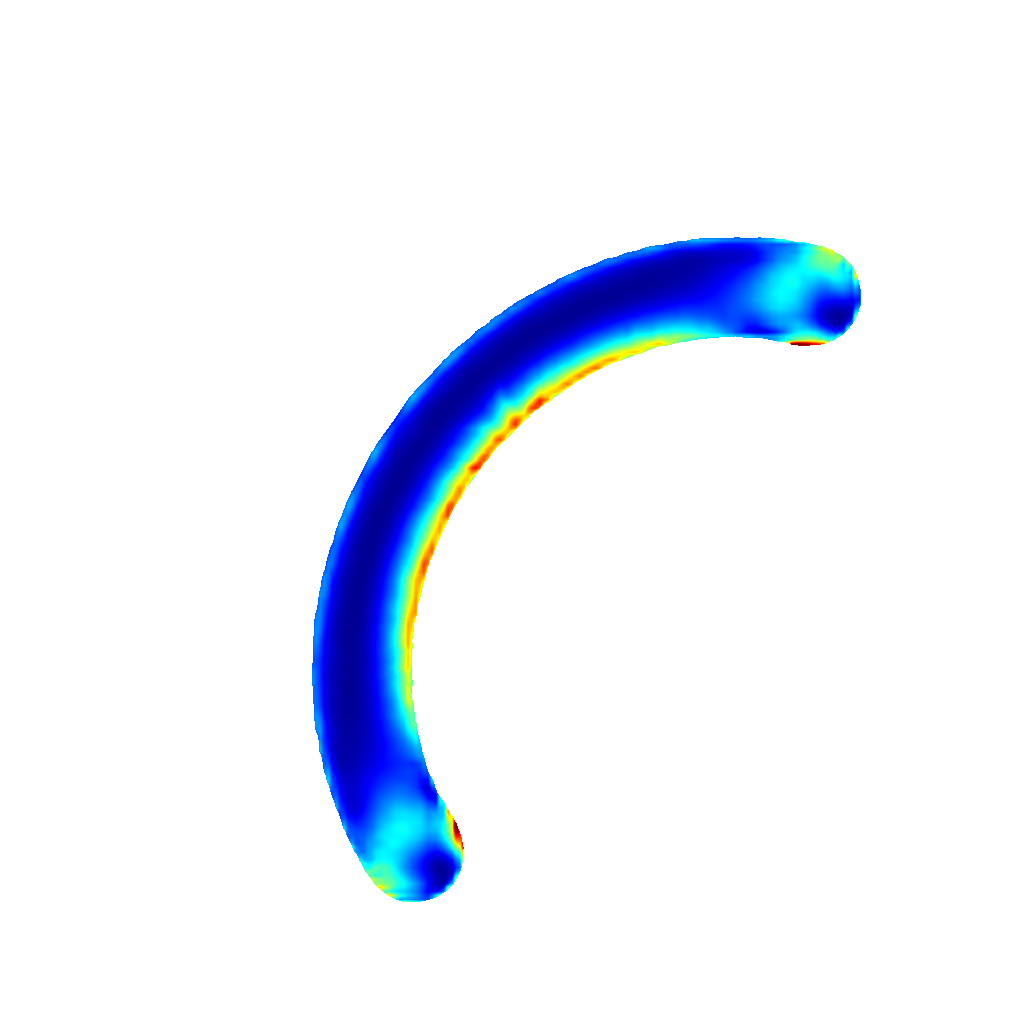}}
\put(.12,-.028){\includegraphics[width=.2\unitlength]{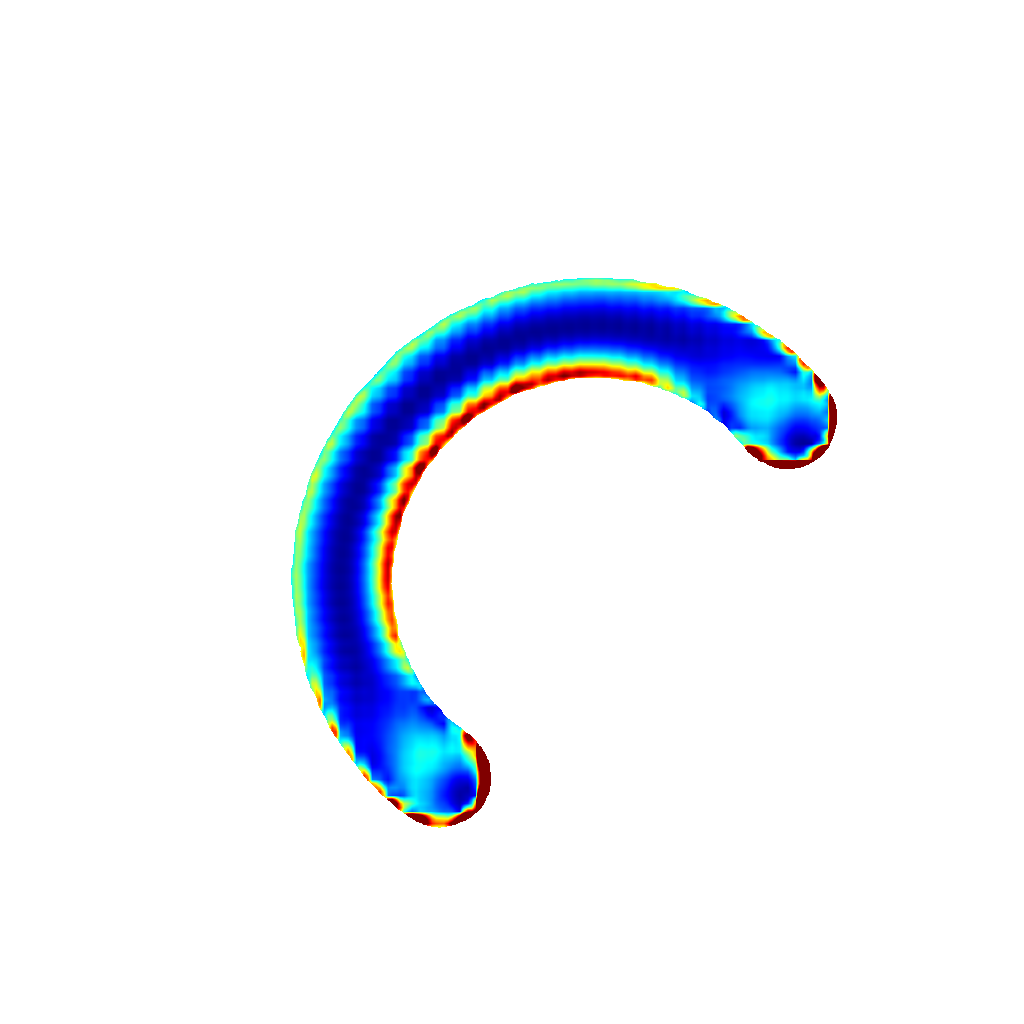}}
\put(.165,-.065){\includegraphics[width=.2\unitlength]{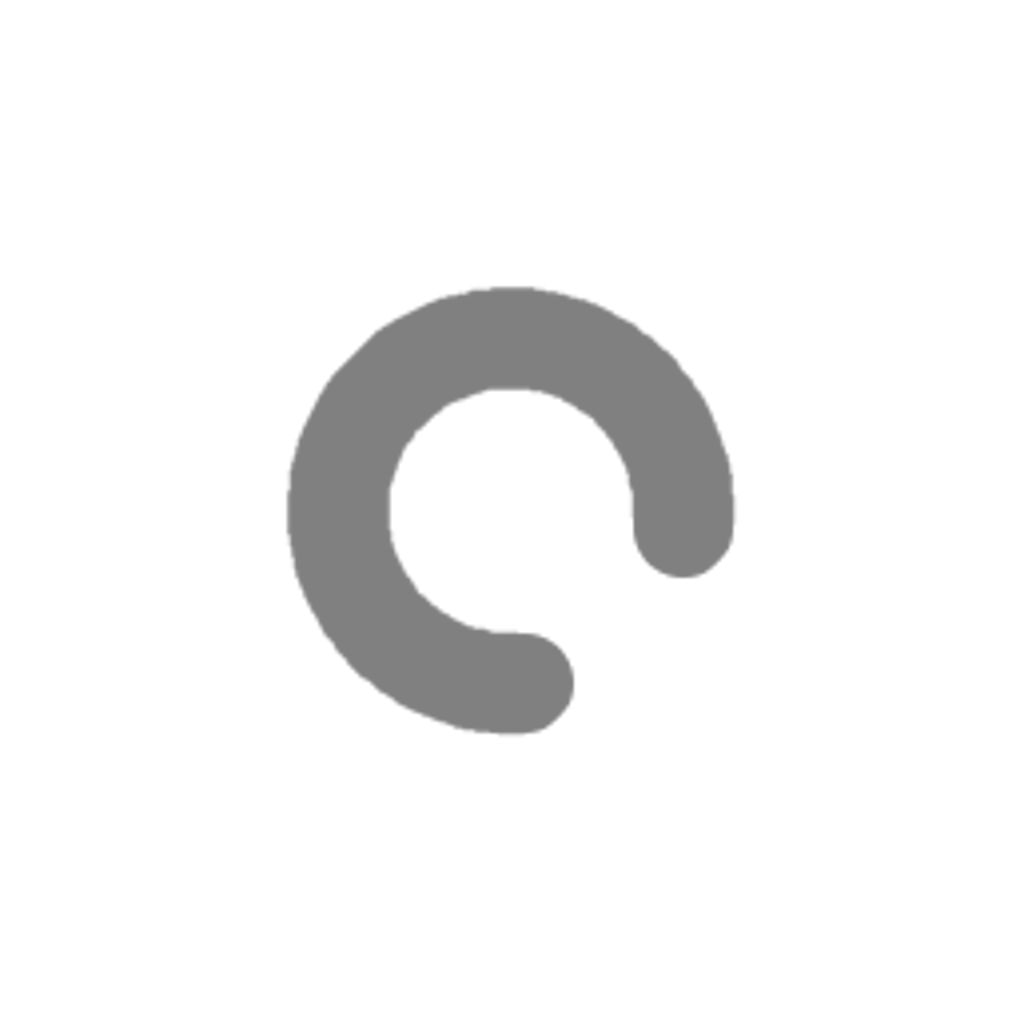}}
\put(.33,.02){\includegraphics[width=.2\unitlength]{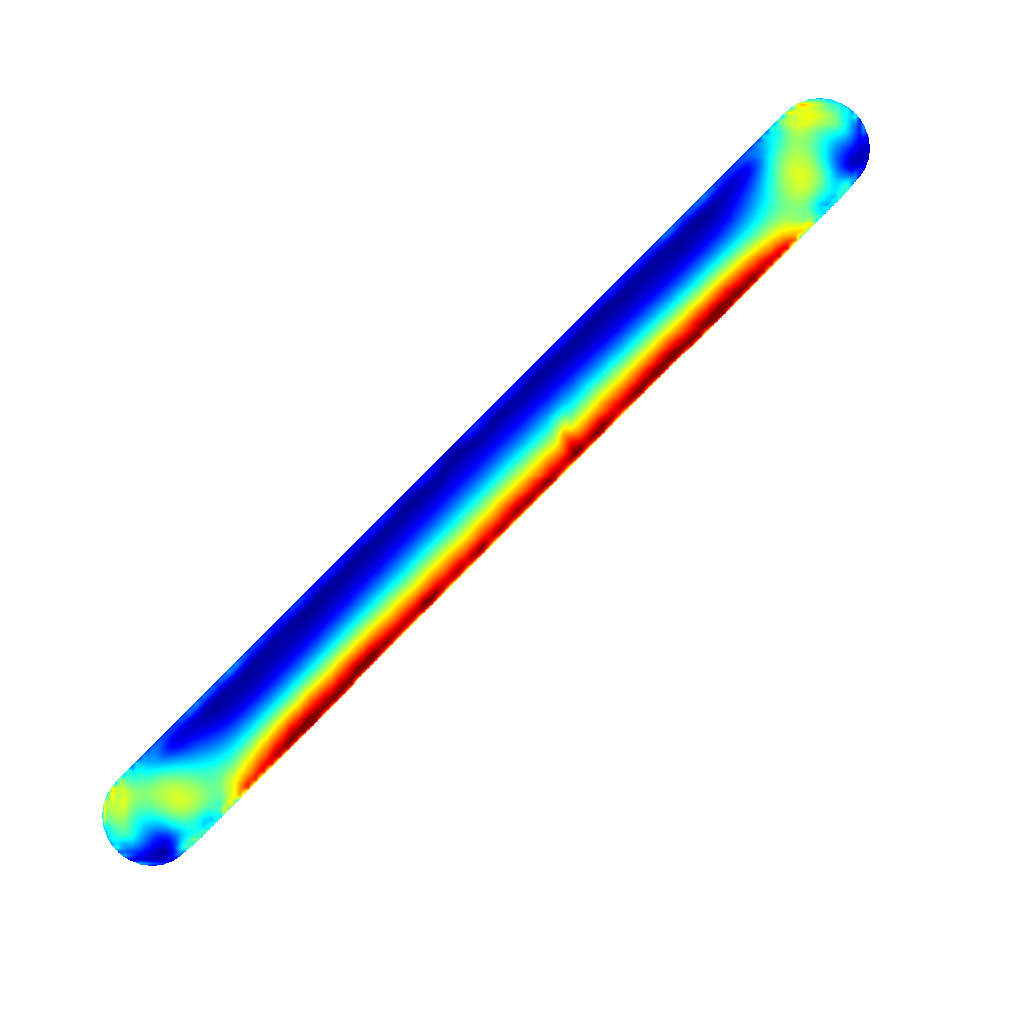}}
\put(.37,.01){\includegraphics[width=.2\unitlength]{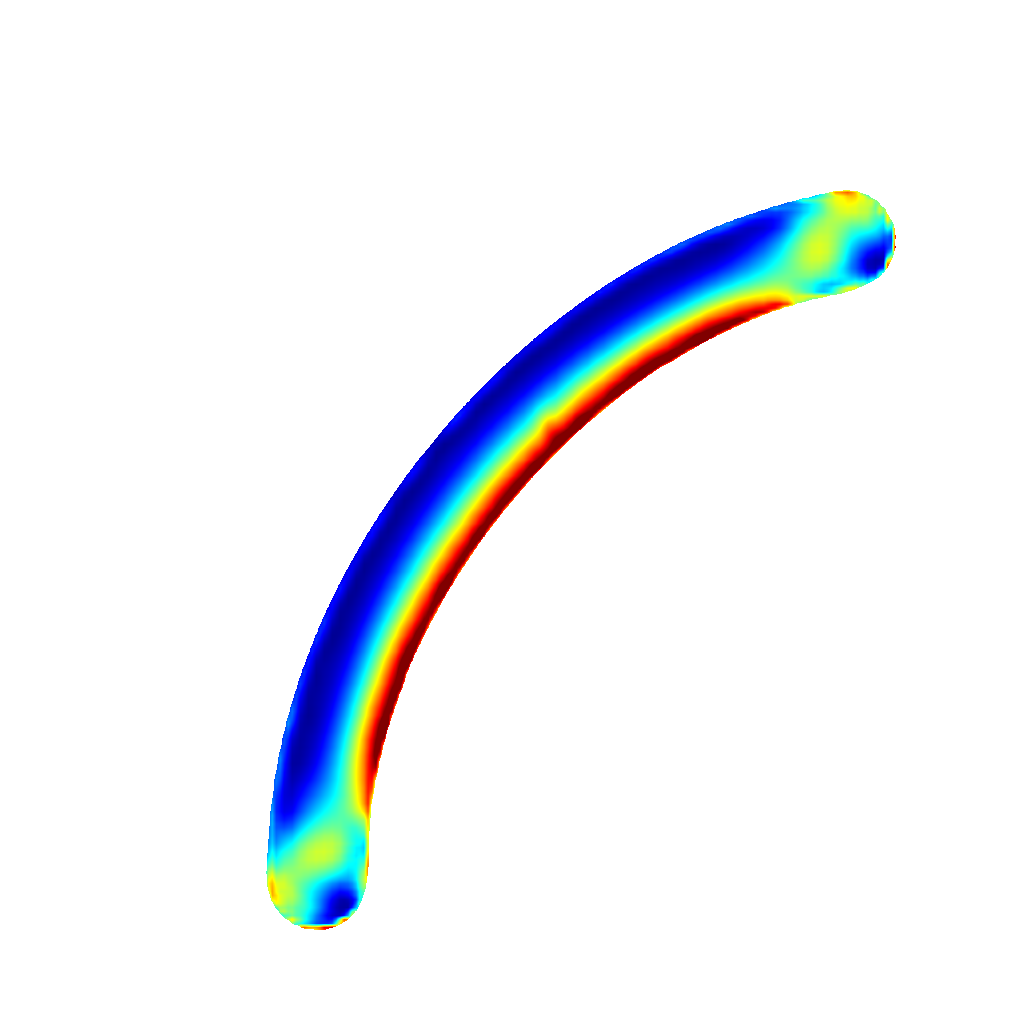}}
\put(.40,-.00){\includegraphics[width=.2\unitlength]{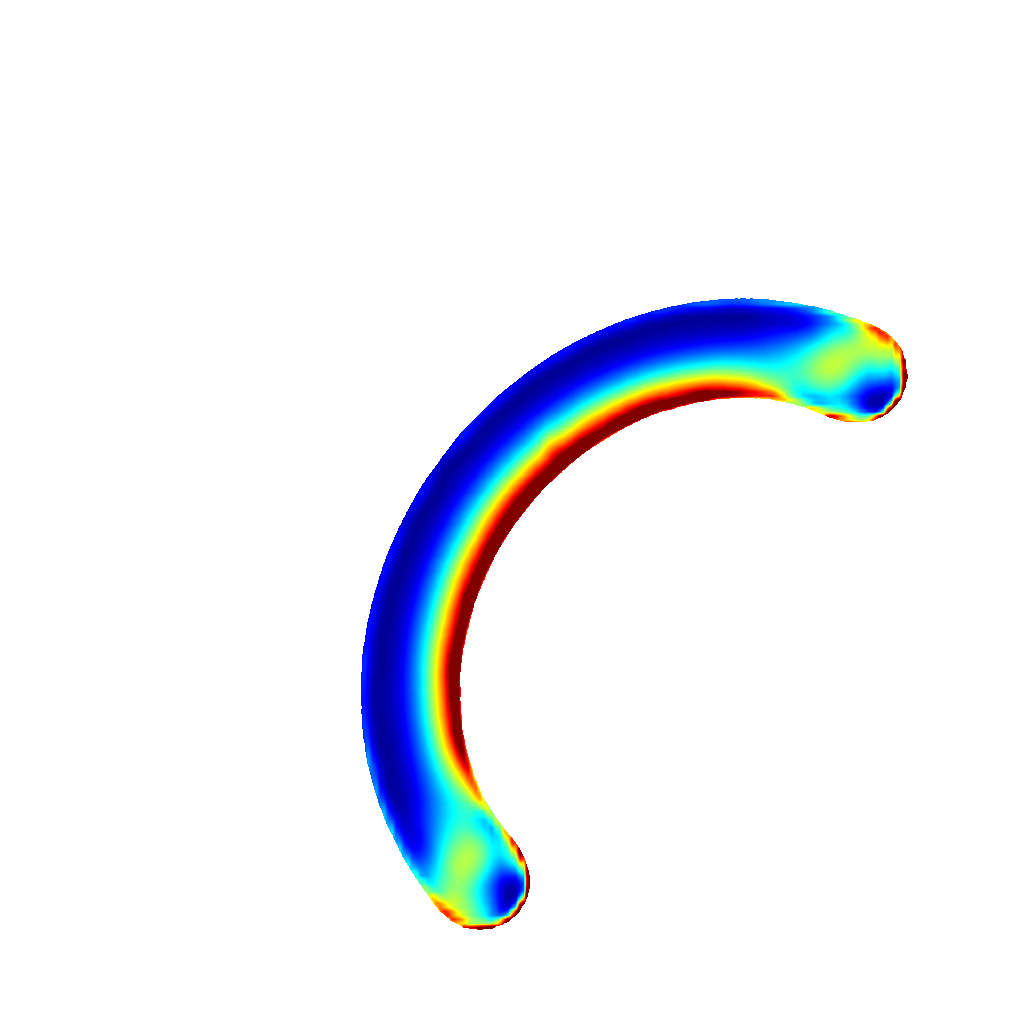}}
\put(.44,-.03){\includegraphics[width=.2\unitlength]{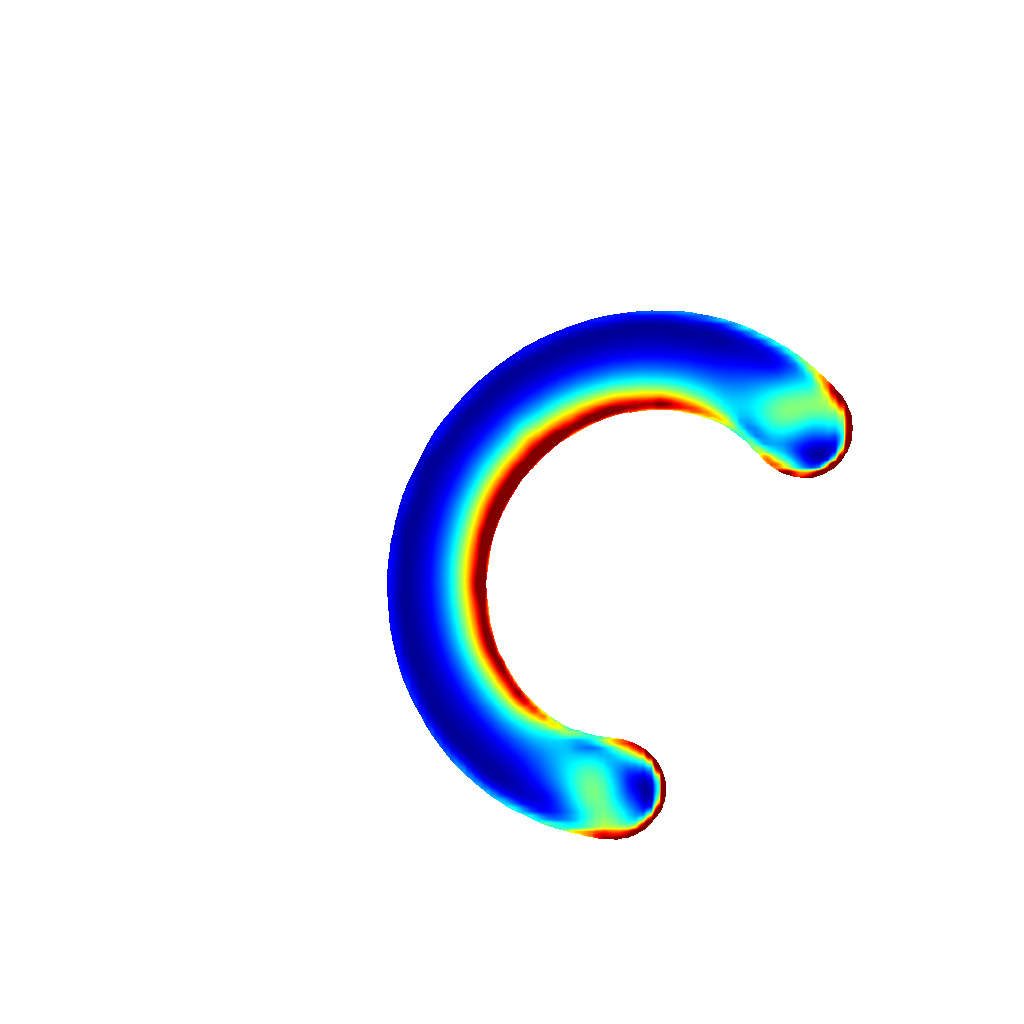}}
\put(.505,-.05){\includegraphics[width=.2\unitlength]{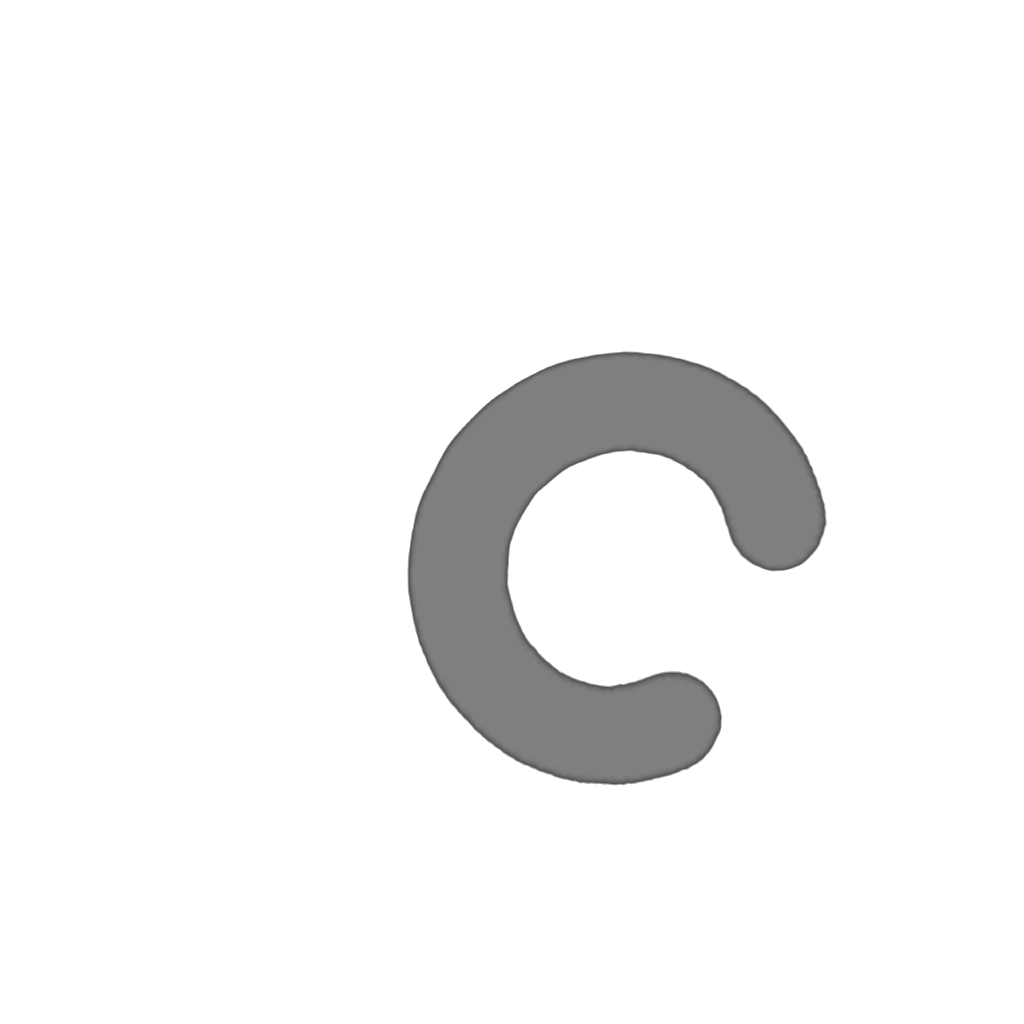}}
\end{picture}
\end{center}
\caption{
Left: Discrete geodesic between two end shapes, where
color indicates the local rate of dissipation (from blue, low, to red, high).
Right: Given the first shape and its variation (in terms of the optimal matching deformation to the second shape), a discrete geodesic is extrapolated.}
\label{fig:shootingGeodesic}
\end{figure}
In fact, discrete geodesics that are variationally described as discrete energy minimizing paths between two given objects can be reproduced via the discrete geodesic flow associated with the discrete exponential map (\conf Figure~\ref{fig:shootingGeodesic}). 

\begin{figure}
\definecolor{lightgrey}{cmyk}{0,0,0,0.1}
\setlength{\unitlength}{.12\linewidth}
\begin{center}
\begin{picture}(7,2)
\put(0,1){\includegraphics[width=.95\unitlength]{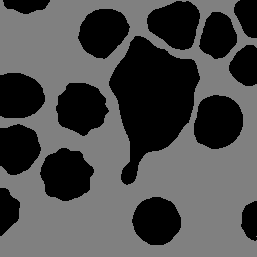}}
\put(1.5,1){\includegraphics[width=.95\unitlength]{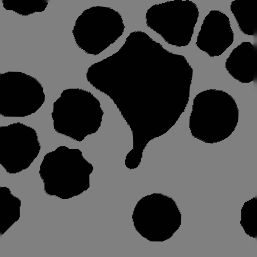}}
\put(0.95,1.4){\Large \color{red} $\boldsymbol{\stackrel{\zeta_1}{\longrightarrow}}$}
\put(3,1){\includegraphics[width=.95\unitlength]{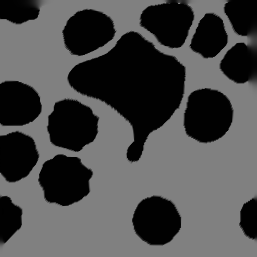}}
\put(4.5,1){\includegraphics[width=.95\unitlength]{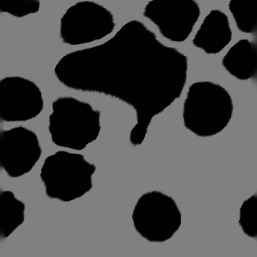}}
\put(6,1){\includegraphics[width=.95\unitlength]{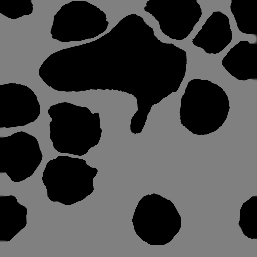}}
\put(0,0){\includegraphics[width=.95\unitlength]{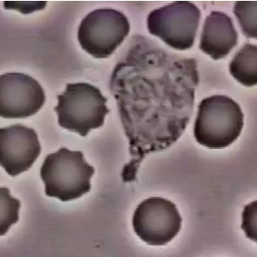}}
\put(1.5,0){\includegraphics[width=.95\unitlength]{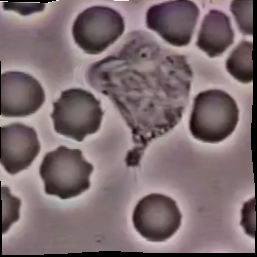}}
\put(3,0){\includegraphics[width=.95\unitlength]{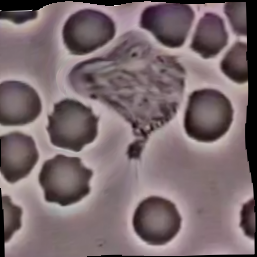}}
\put(4.5,0){\includegraphics[width=.95\unitlength]{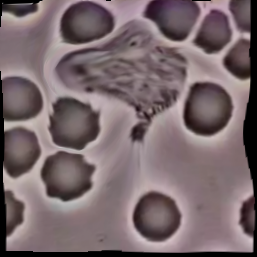}}
\put(6,0){\includegraphics[width=.95\unitlength]{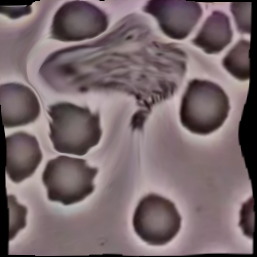}}
\end{picture}
\end{center}
\caption{\label{fig:bloodextrapolate} Top: Given the first shape and its initial variation $\zeta_1$ represented by the second shape, a discrete geodesic is extrapolated.
Bottom: The texture of a video frame can be transported along with the varying shapes.}
\end{figure}

As for the discrete logarithm we experimentally observe convergence of the discrete exponential map  in the sense
\begin{equation}\label{eq:EXPconvergence}
\ExpO{k}{\object} \left(\textstyle \frac1k \zeta\right) \to \exp_\object(\zeta)\quad\mbox{for $k\to\infty$}
\end{equation}
as shown in Figure~\ref{fig:LogExp}. An example of geodesic shape extrapolation for multicomponent objects is depicted in Figure~\ref{fig:bloodextrapolate}.

\subsection{Discrete parallel transport and detail transfer}
\begin{figure}
\begin{center}
\setlength{\unitlength}{0.6\linewidth}
\begin{picture}(1,.45)
\put(0,.05){\includegraphics[width=.9\unitlength,angle=-20,origin=c]{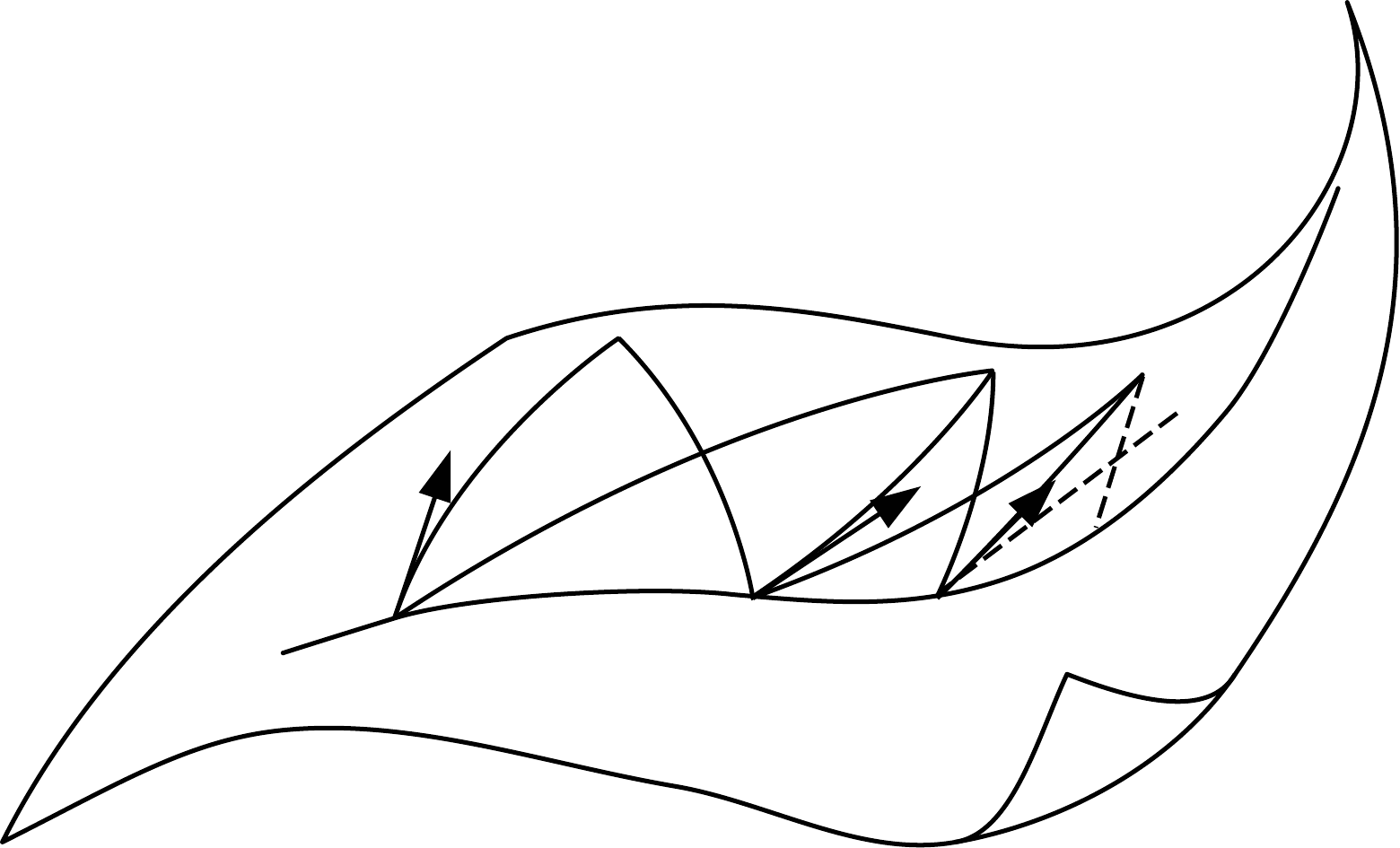}}
\put(.1,.19){$\manifold$}
\put(.8,.16){$\curve(t)$}
\put(.23,.21){$\object(\!(k\!-\!1)\tau)$}
\put(.45,.14){$\object(k\tau)$}
\put(.49,.415){$\object^p_{k-1}$}
\put(.71,.31){$\object^p_k$}
\put(.515,.33){$\object^{\times}_k$}
\put(.499,.295){$\bullet$}
\put(.24,.33){$\zeta_{k-1}$}
\put(.55,.24){$\zeta_{k}$}
\put(.28,.26){$\bullet$}
\put(.475,.38){$\bullet$}
\put(.5,.19){$\bullet$}
\put(.695,.28){$\bullet$}
\end{picture}
\end{center}
\caption{
A sketch of the  discrete parallel transport of $\zeta_{k-1}\in T_{\object((k-1)\tau)}\manifold$ via Schild's ladder along the edge from $\object((k-1)\tau)$ to $\object(k\tau)$ on a  curve in shape space.}
\label{fig:SchildLadder}
\end{figure}

\begin{figure}
\setlength{\unitlength}{.09\linewidth}
\begin{picture}(5,1.85)(.3,.25)
\put(0,.99){\includegraphics[width=1.5\unitlength]{PA1}}
\put(1,.96){\includegraphics[width=1.5\unitlength,angle=-8,origin=c]{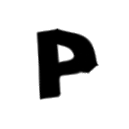}}
\put(2,.96){\includegraphics[width=1.5\unitlength,angle=-13,origin=c]{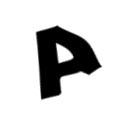}}
\put(2.9,.92){\includegraphics[width=1.5\unitlength,angle=-23,origin=c]{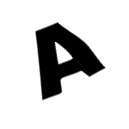}}
\put(4.2,.99){\includegraphics[width=1.5\unitlength]{PA5}}
\put(0,0){\includegraphics[width=1.5\unitlength]{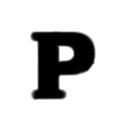}}
\put(1,-.03){\includegraphics[width=1.5\unitlength,angle=-5,origin=c]{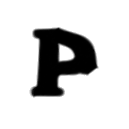}}
\put(2,-.06){\includegraphics[width=1.5\unitlength,angle=-10,origin=c]{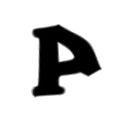}}
\put(3,-.08){\includegraphics[width=1.5\unitlength,angle=-15,origin=c]{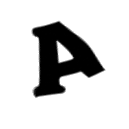}}
\put(4,-.1){\includegraphics[width=1.5\unitlength,angle=-20,origin=c]{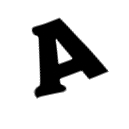}}
\put(.41,.38){\color{red}\circle{.18}}
\put(.41,1.07){\color{red}\circle{.18}}
\put(.715,.38){\color{red}\circle{.18}}
\end{picture}
\hspace{3ex}
\setlength{\unitlength}{.1\linewidth}
\begin{picture}(5,1.9)(0,.1)
\put(0,1){\includegraphics[width=1.2\unitlength,angle=0,origin=c]{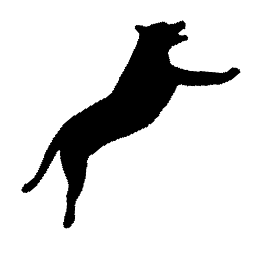}}
\put(1,.98){\includegraphics[width=1.2\unitlength,angle=0,origin=c]{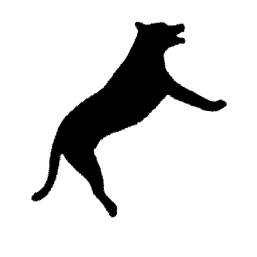}}
\put(2,.96){\includegraphics[width=1.2\unitlength,angle=0,origin=c]{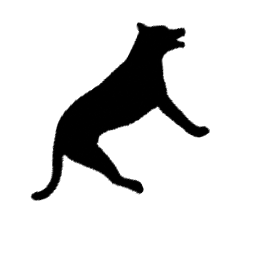}}
\put(3,1.02){\includegraphics[width=1.2\unitlength,angle=0,origin=c]{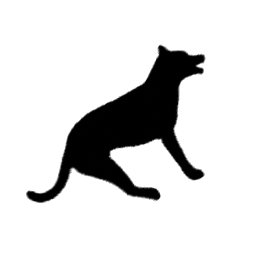}}
\put(4,1){\includegraphics[width=1.2\unitlength,angle=0,origin=c]{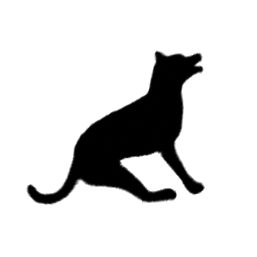}}
\put(0,0){\includegraphics[width=1.2\unitlength]{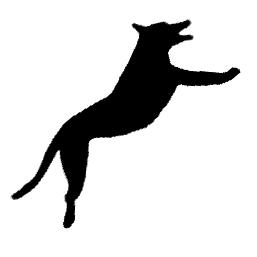}}
\put(1,-.01){\includegraphics[width=1.2\unitlength,angle=0,origin=c]{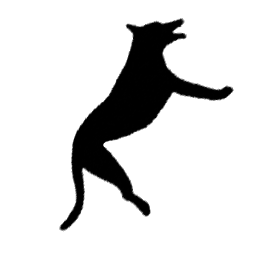}}
\put(2,-.02){\includegraphics[width=1.2\unitlength,angle=0,origin=c]{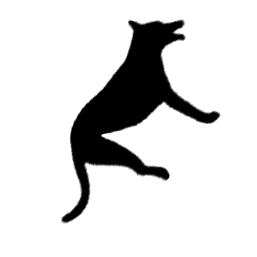}}
\put(3,-.03){\includegraphics[width=1.2\unitlength,angle=0,origin=c]{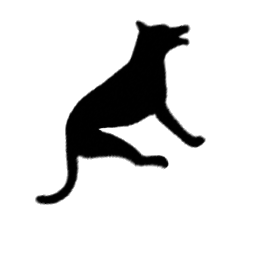}}
\put(4,0){\includegraphics[width=1.2\unitlength,angle=0,origin=c]{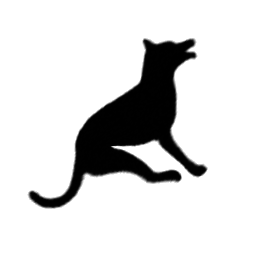}}
\put(.08,.28){\color{red}\circle{.15}}
\put(.42,.5){\color{red}\circle{.15}}
\put(.34,.23){\color{red}\circle{.15}}
\put(.67,1.07){\color{red}\circle{.15}}
\put(.9,1.05){\color{red}\circle{.15}}
\end{picture}
\caption{Discrete parallel transport is applied along two discrete geodesic paths 
connecting letters 'P' and 'A' (left) and two different poses of a dog (right), respectively.
On the top the original discrete geodesics are shown, while the bottom left shows the discrete
parallel transport of serifs along the geodesic between the two letters, and the bottom right shows the transport of changes on the first dog's shape, which allows to copy the changes to the other poses.}
\label{fig:parallelTransport}
\end{figure}

Parallel transport allows to translate a vector $\zeta \in T_\object \manifold$ (which is considered as the variation of an object $\object=\object(0)$)
along a curve $(\object(t))_{t \in [0,1]}$ in shape space. The resulting
$(\zeta(t))_{t \in [0,1]}$ changes as little as possible while keeping the angle between $\zeta(t)$ and the path velocity $v(t)$ fixed. Using the Levi-Civita connection this can be phrased as $\nabla_{v(t)} \zeta(t)=0$.
There is a well-known first-order approximation of parallel transport called Schild's ladder \cite{EhPiSc72,KhMiNe00},
which is based on the construction of a sequence of geodesic parallelograms, sketched in Figure~\ref{fig:SchildLadder}, where the two diagonal geodesics always meet at their midpoints.
Given a curve $(\object(t))_{t \in [0,1]}$ and a tangent vector $\zeta_{k-1}\in T_{\object((k-1)\tau)}\manifold$,
the approximation $\zeta_k\in T_{\object(k \tau)}\manifold$ of the parallel transported vector via a geodesic parallelogram can be expressed as
\begin{align*}
\object^p_{k-1} &= \exp_{\object((k-1)\tau)} \zeta_{k-1} \,,\\
\object_k^{\times} &= \exp_{\object^p_{k-1}} \frac12 \log_{\object^p_{k-1}}  \object(k \tau) \,,\\
\object_{k}^p &= \exp_{\object((k-1)\tau)} 2 \log_{\object((k-1)\tau)}  \object_k^{\times} \,,\\
\zeta_k &= \log_{\object(k \tau)} \object_{k}^p\,.
\end{align*}
Here, $\object_k^{\times}$ is the midpoint of the two diagonals of the geodesic parallogramm with vertices $\object((k-1)\tau)$,
$\object^p_{k-1}$, $\object^p_k$, and $\object(k\tau)$.
This scheme can be easily transferred to discrete curves $(\object_0,\ldots, \object_K)$ in shape space based on the discrete logarithm and the discrete exponential introduced above.
In the $k$th step of the discrete transport we start with a displacement $\zeta_{k-1}$ on
$\object_{k-1}$ and compute
\begin{align*}
\object^p_{k-1} &= \ExpO{1}{\object_{k-1}} \zeta_{k-1}\,, \\
\object_k^{\times} &= \ExpO{1}{\object^p_{k-1}} \left( \Log{2}_{\object^p_{k-1}}  (\object_{k})\right)\,, \\
\object_k^p &= \ExpO{2}{\object_{k-1}} \left(\Log{1}_{\object_{k-1}}  (\object_k^{\times})  \right)\,, \\
\zeta_{k} &= \Log{1}_{\object_k} (\object_k^p)\,,
\end{align*}
where $\zeta_{k}$ is the transported displacement on $\object_k$. 
Here, $\object_k^{\times}$ is the midpoint of the two discrete geodesics of order $2$ with end points $\object^p_{k-1}$,
$\object_k$ and $\object_{k-1}$, $\object_k^p$, respectively.
Since the last of the above steps is the inverse of the first step in the subsequent iteration, these steps need to be performed only for $k=K$. We will denote the resulting transport operator by $P_{\object_K,\ldots, \object_0}$.
Figure~\ref{fig:parallelTransport} shows examples of discrete parallel transport for feature transfer along curves in shape space.

{\bf Remark:} As in the continuous case, the discrete parallel transport can be used to define a discrete Levi-Civita connection. For $\xi \in T_\object \manifold$ and for a vector field $\eta$ in the tangent bundle $T\manifold$ one computes $\object_\tau= \ExpO{1}{\object}(\tau \xi)$ and then defines
\notinclude{
$$\nabla_{\xi}^\tau \eta := \frac1\tau \left(P_{\object_\tau, \object}^{-1} \eta(\object_\tau)-\eta(\object)\right)$$
as the time discrete connection with time step size $\tau$,
where this particular form was chosen so as to make the displacements $\zeta_k$ along a discrete geodesic satisfy $\nabla_{\zeta_k}^1\zeta_k=0$.
}
$$\nabla_{\xi}^\tau \eta := \frac1\tau \left(P_{\object, \object_\tau} \eta(\object_\tau)-\eta(\object)\right)$$
as the time discrete connection with time step size $\tau$.
\notinclude{, where we have used that $P_{\object_\tau, \object}^{-1}$
equals $P_{\object, \object_\tau}$ up to higher order terms in $\tau$.
}

\section{Numerical discretization}
\label{sec:numerics}
The proposed discrete geodesic calculus requires an effective and efficient spatial discretization of 
\begin{itemize}
\item[-]
volumetric objects $\object$ in the underlying shape space, 
\item[-]
of nonlinear deformations $\psi$ to encode matching correspondences, 
\item[-]
and of linear displacements $\zeta$ as approximate tangent vectors.  
\end{itemize}
We restrict ourselves here to the case of objects $\object \subset \R^2$. To this end we consider the space $\V_h$ of piecewise affine finite element functions on a regular simplicial mesh over a rectangular computational domain $\compDom$. Here $h$ indicates the grid size, where in our applications $h$ ranges from a coarse grid size $2^{-6}$ to a fine grid size $2^{-8}$.
Then, deformations and displacements are considered as functions in $(\V_h)^2$.
Objects $\object$, the original degrees of freedom in our geometric calculus,
will be represented via deformations $\phi$ over reference objects $\hat \object$  (\eg $\refobject$), \ie 
$\object=\phi(\hat\object)$. These reference objects are encoded by approximate characteristic functions $\chi^h_{\hat\object}\in\V_h$ and the deformations $\phi$ are considered as injective deformations $\phi: \compDom \to \R^2$ and discretized as elements in $(\V_h)^2$.

\subsection{Parameterization of discrete geodesics}
\label{sec:parageodesic}
\begin{figure}
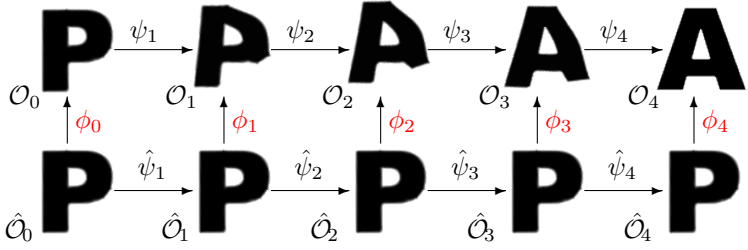

\begin{center}
\setlength{\unitlength}{.08\linewidth}
\begin{picture}(10,2.9)(.4,.4)
\put(0.5,0){\includegraphics[width=2\unitlength]{PA1}}
\put(2.5,0){\includegraphics[width=2\unitlength]{PA1}}
\put(4.5,0){\includegraphics[width=2\unitlength]{PA1}}
\put(6.5,0){\includegraphics[width=2\unitlength]{PA1}}
\put(8.5,0){\includegraphics[width=2\unitlength]{PA1}}
\put(0.6,.4){$\hat \object_{\!0}$}
\put(2.6,.4){$\hat \object_{\!1}$}
\put(4.5,.4){$\hat \object_{\!2}$}
\put(6.5,.4){$\hat \object_{\!3}$}
\put(8.5,.4){$\hat \object_{\!4}$}
\put(2,1){\vector(1,0){1}}
\put(4,1){\vector(1,0){1}}
\put(6,1){\vector(1,0){1}}
\put(8,1){\vector(1,0){1}}
\put(2.3,1.2){$\hat\psi_1$}
\put(4.3,1.2){$\hat\psi_2$}
\put(6.3,1.2){$\hat\psi_3$}
\put(8.3,1.2){$\hat\psi_4$}
\put(1.4,1.6){\vector(0,1){.6}}
\put(3.4,1.6){\vector(0,1){.6}}
\put(5.4,1.6){\vector(0,1){.6}}
\put(7.4,1.6){\vector(0,1){.6}}
\put(9.4,1.6){\vector(0,1){.6}}
\put(1.5,1.8){\textcolor{red}{$\phi_0$}}
\put(3.5,1.8){\textcolor{red}{$\phi_1$}}
\put(5.5,1.8){\textcolor{red}{$\phi_2$}}
\put(7.5,1.8){\textcolor{red}{$\phi_3$}}
\put(9.5,1.8){\textcolor{red}{$\phi_4$}}
\put(0.50,1.8){\includegraphics[width=2\unitlength]{PA1}}
\put(2.35,1.8){\includegraphics[width=2\unitlength,angle=-8,origin=c]{PA2}}
\put(4.25,1.8){\includegraphics[width=2\unitlength,angle=-15,origin=c]{PA3}}
\put(6.20,1.7){\includegraphics[width=2\unitlength,angle=-22,origin=c]{PA4}}
\put(8.50,1.8){\includegraphics[width=2\unitlength]{PA5}}
\put(0.65,2.1){$\object_0$}
\put(2.65,2.1){$\object_1$}
\put(4.65,2.1){$\object_2$}
\put(6.65,2.1){$\object_3$}
\put(8.55,2.1){$\object_4$}
\put(2,2.8){\vector(1,0){1}}
\put(4,2.8){\vector(1,0){1}}
\put(6,2.8){\vector(1,0){1}}
\put(8,2.8){\vector(1,0){1}}
\put(2.2,2.95){$\psi_1$}
\put(4.2,2.95){$\psi_2$}
\put(6.2,2.95){$\psi_3$}
\put(8.2,2.95){$\psi_4$}
\end{picture}
\end{center}
\caption{A diagram illustrating the parameterization of the domains $\object_k$ along the discrete geodesic via deformations $\phi_k$ --- indicated in red as the actual degrees of freedom --- over reference domains $\hat \object_k$.}
\label{fig:commutativeDiagram}
\end{figure}

To compute a discrete geodesic --- different from \cite{WiBaRu10} --- we now replace the objects $\object_0, \ldots, \object_K$ as arguments of the energy \eqref{eqn:discretePathEnergy} by associated deformations $\phi_0, \ldots, \phi_K$
over a set of reference domains $\hat \object_0, \ldots, \hat \object_K$ as described above.
By this technique, instead of $K$ deformations and $K-1$ domain descriptions (\eg via level sets) as in \cite{WiBaRu10} we will be able to consider solely $K+1$ parameterizing deformations, which turns out to be a significant computational advantage.
Next, we assume that reference matching deformations $\hat\psi_1,\ldots, \hat\psi_K$
are given with $\hat \object_{k}=\hat \psi_k(\hat \object_{k-1})$
(\conf Figure~\ref{fig:commutativeDiagram}). 
Now, we express the matching deformations $(\psi_k)_{k=1,\ldots, K}$ over which we minimize in \eqref{eqn:discretePathEnergy} in terms of the parameterizing deformations $(\phi_k)_{k=0,\ldots, K}$ and the reference matching deformations $(\hat \psi_k)_{k=1,\ldots, K}$ and set
$$\psi_k= \phi_k \circ \hat \psi_k \circ \phi_{k-1}^{-1}$$
for $k=1,\ldots, K$. Now, using a change of variables one can rewrite the deformation energies $\W_{\object_\ikm}[\psi_\ik]$ in
\eqref{eqn:discretePathEnergy} as
\begin{equation*}
\int_D \chi_{\hat \object_{\ikm}} W(\nabla(\phi_\ik \circ \hat \psi_\ik)(\nabla \phi_\ikm)^{-1}) \det \nabla \phi_\ikm \d \hat x\,.
\end{equation*}
Furthermore, instead of void we consider a $\delta_1$ times softer material outside the object domain replacing $\chi_{\hat\object}$ by $\chi^{\delta_1}_{\hat\object} =  (1-\delta_1) \chi_{\hat\object} + \delta_1$
so that altogether the deformation energy $\W_{\object_\ikm}[\psi_\ik]$ is replaced by the following energy over the parameterizing deformations $\phi_\ikm$ and $\phi_\ik$:
\begin{equation*}
\W_{\hat\object_\ikm}^{\delta_1,\hat\psi_\ik}[\phi_\ikm,\phi_\ik]=
\int_D\chi^{\delta_1}_{\hat\object_{\ikm}} W(\nabla(\phi_\ik \circ \hat \psi_\ik)(\nabla \phi_\ikm)^{-1}) \det \nabla\phi_\ikm \d \hat x
\end{equation*}
The condition that $\object_0$ and $\object_K$ are prescribed is taken care of with
penalty functionals
\begin{equation*}
\F^{\varepsilon,\delta_2}_{\hat\object_i,\object_i}[\phi_i]= \frac1\varepsilon
\int_D(G_{\delta_2}*\chi_{\hat\object_i}-G_{\delta_2}*\chi_{\object_i}\circ \phi_i)^2\,\d \hat x
\end{equation*}
for $i=0,K$, where $G_{\delta_2}$ is a Gaussian of filter width $\delta_2$.
Finally, to ensure that not only the concatenations $\phi_k \circ \hat \psi_k \circ \phi_{k-1}^{-1}$ of deformations are regular but also every single deformation $\phi_k$, we add a term $\delta_3 \W_D[\phi_k]$ for all $k=0, \ldots, K$.
Summarizing, to compute a discrete geodesic between two shapes $\object_0,\object_K$ or a discrete logarithm,
we minimize the total energy
\begin{equation*}
\sum_{\ik=1}^K\W_{\hat\object_\ikm}^{\delta_1,\hat\psi_\ik}[\phi_\ikm,\phi_\ik]
+\delta_3\sum_{\ik=0}^K\W_D[\phi_\ik]
+\F^{\varepsilon,\delta_2}_{\hat\object_0,\object_0}[\phi_0]+\F^{\varepsilon,\delta_2}_{\hat\object_K,\object_K}[\phi_K]
\end{equation*}
over all the parameterizing deformations $\phi_0,\ldots,\phi_K$. This minimization then determines the shapes $\object_k=\phi_k(\hat\object_k)$ forming the discrete geodesic
as well as the discrete logarithm 
$$\Log{K}_{\object_0}(\object_K)=\phi_1 \circ \hat \psi_1 \circ \phi_{0}^{-1}-\id\,.$$
Due to the assumptions on the energy integrand $W$, self-penetration of deformed objects is ruled out in our fully discrete model, which as a byproduct leads to topology preservation along discrete geodesics.
The resulting matching deformations $\psi_1,\ldots,\psi_K$ are determined up to translations and rotations due to the built-in frame indifference (\conf Section\,\ref{sec:approxdistance}).
In our applications the parameters of the algorithm are chosen as follows:
$\delta_1=0.01$, $\delta_2=h$, $\delta_3=0.01$, $\varepsilon=0.1$.

Furthermore, the previously described finite element discretization is applied, and
the energies are computed via Simpson quadrature on each element.
As opposed to Gaussian quadrature, this quadrature rule has the advantage that quadrature points lie in the corners of each finite element,
which is where the extremal values of the finite element functions and their gradients occur.
This is relatively important when dealing with deformations as finite element functions,
since a Gaussian quadrature rule might for example miss that a deformation exhibits self-penetration in the corner of one element.
Such deformations will thus not be rejected during energy minimization,
which in turn results in technical difficulties when computing pullbacks or pushforwards with respect to these deformations
and when prolongating them onto grids with finer resolution.

Pullbacks of a finite element function $f\in\V_h$ with respect to a finite element deformation $\phi\in(\V_h)^2$ at quadrature points $x$ are computed
by first evaluating the deformation at that point.
If $\phi(x)$ lies outside the computational domain $D$, it is projected back onto the boundary $\partial D$
(this approach is more robust than \eg neglecting such points,
since it avoids structural changes when deformed quadrature points toggle between inside and outside).
Finally, $f$ is evaluated at the projected position.

The minimization is performed by a Newton trust region algorithm \cite{Conn2000} on the $(K+1)$-tuple of deformations
$(\phi_0,\ldots, \phi_K)$ and requires the evaluation of first and second derivatives of the energy.
Note that the second derivative involves mixed derivatives with respect to $\phi_{k-1}$ and $\phi_k$ for all $k=1,\ldots, K$.
During the Newton iteration, these mixed terms provide a coupling between all deformations
which results in a fast relaxation and balance between them.

\begin{figure}
\begin{center}
\setlength{\unitlength}{.08\linewidth}
\begin{picture}(11,3)(.4,0)
\put(.4,0.4){$K\!=\!8$}
\put(.4,1.4){$K\!=\!4$}
\put(.4,2.4){$K\!=\!2$}
\put(2.5,2.5){\line(1,0){3}}
\put(6.5,2.5){\line(1,0){3}}
\put(2.5,1.5){\line(1,0){1}}
\put(4.5,1.5){\line(1,0){1}}
\put(6.5,1.5){\line(1,0){1}}
\put(8.5,1.5){\line(1,0){1}}
\put(2.2,0.5){\line(1,0){0.4}}
\put(3.2,0.5){\line(1,0){0.4}}
\put(4.2,0.5){\line(1,0){0.4}}
\put(5.2,0.5){\line(1,0){0.4}}
\put(6.3,0.5){\line(1,0){0.4}}
\put(7.3,0.5){\line(1,0){0.4}}
\put(8.4,0.5){\line(1,0){0.4}}
\put(9.35,0.5){\line(1,0){0.4}}
\put(1.5,0){\includegraphics[width=\unitlength]{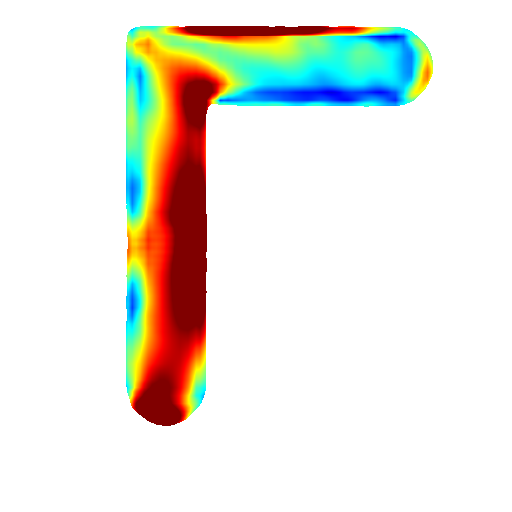}}
\put(2.5,0){\includegraphics[width=\unitlength]{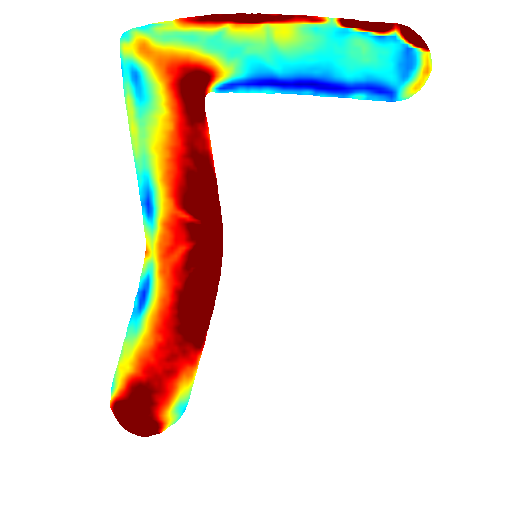}}
\put(3.5,0){\includegraphics[width=\unitlength]{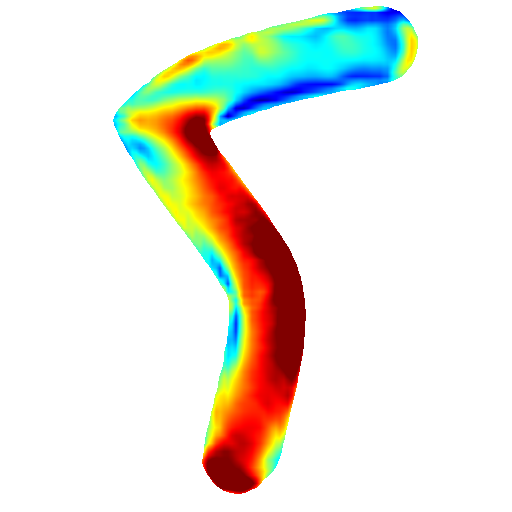}}
\put(4.5,0){\includegraphics[width=\unitlength]{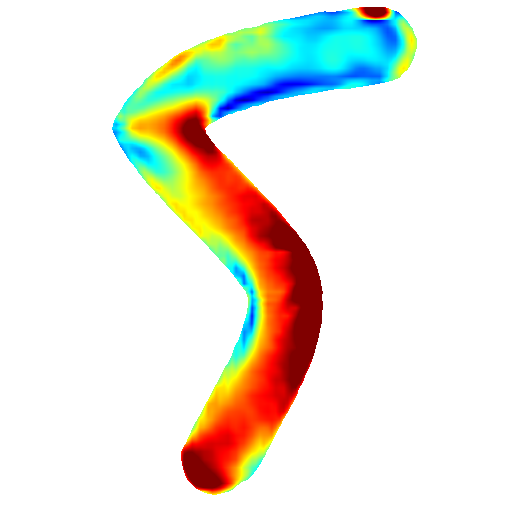}}
\put(5.5,0){\includegraphics[width=\unitlength]{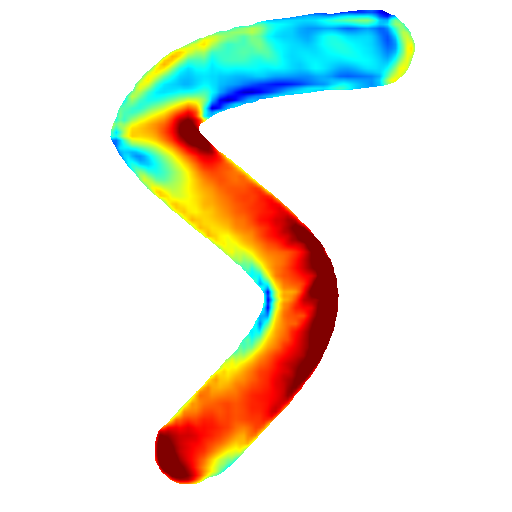}}
\put(6.5,0){\includegraphics[width=\unitlength]{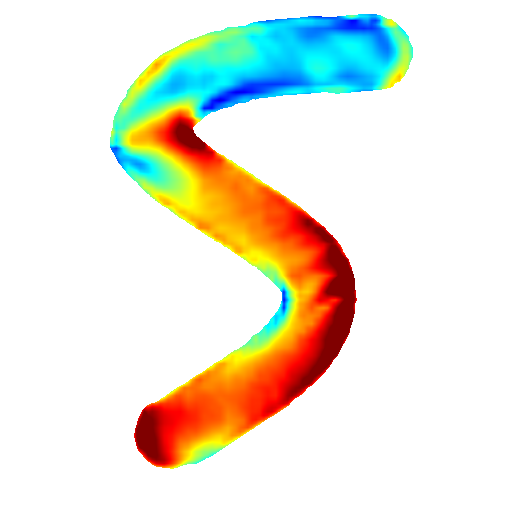}}
\put(7.5,0){\includegraphics[width=\unitlength]{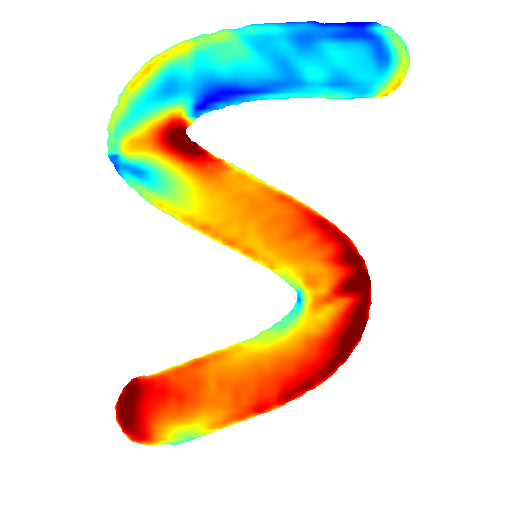}}
\put(8.5,0){\includegraphics[width=\unitlength]{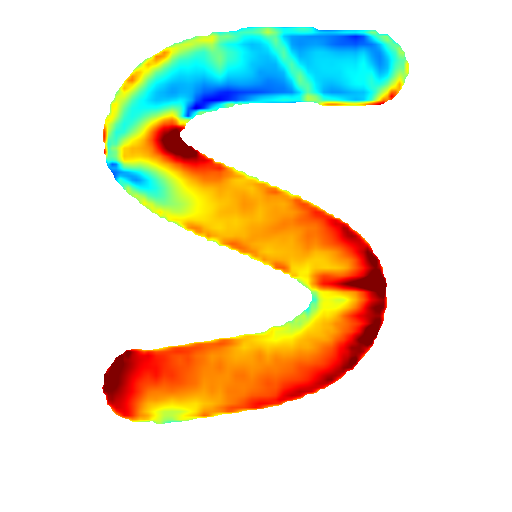}}
\put(9.5,0){\includegraphics[width=\unitlength]{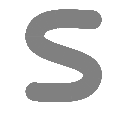}}
\put(1.5,1){\includegraphics[width=\unitlength]{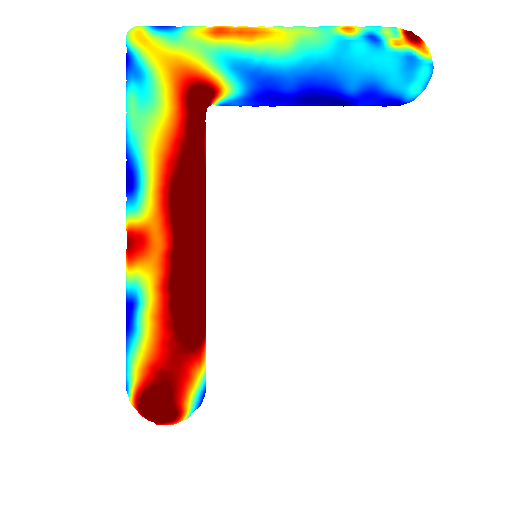}}
\put(3.5,1){\includegraphics[width=\unitlength]{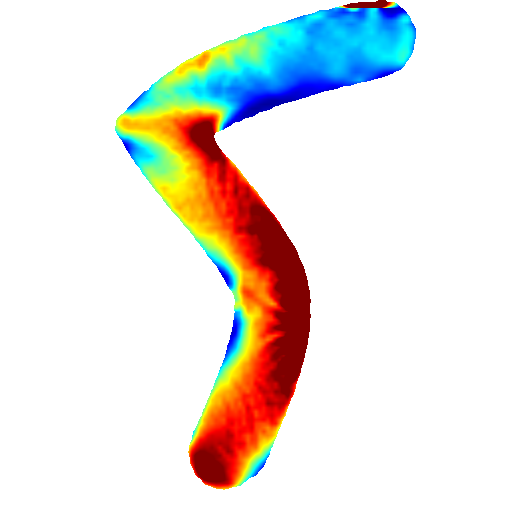}}
\put(5.5,1){\includegraphics[width=\unitlength]{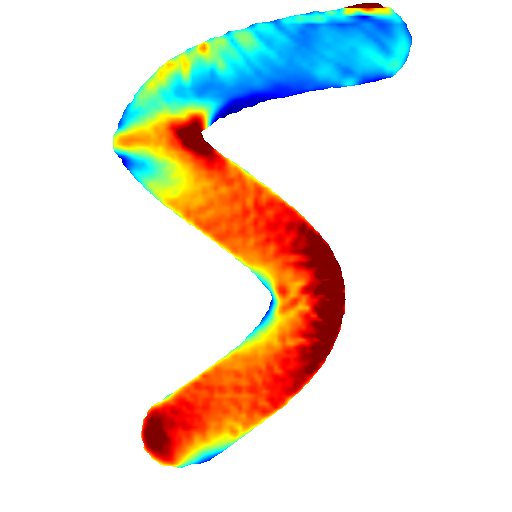}}
\put(7.5,1){\includegraphics[width=\unitlength]{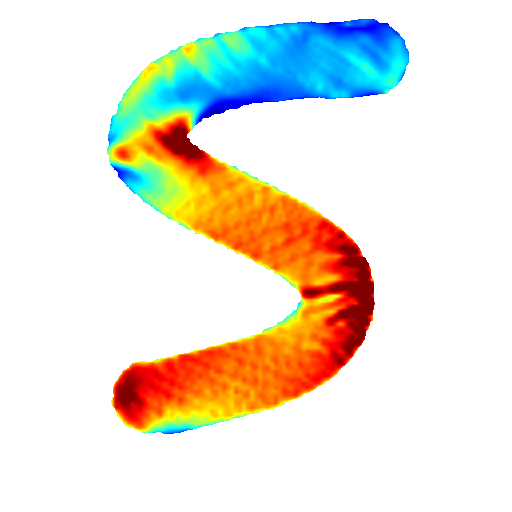}}
\put(9.5,1){\includegraphics[width=\unitlength]{GammaSCol_S}}
\put(1.5,2){\includegraphics[width=\unitlength]{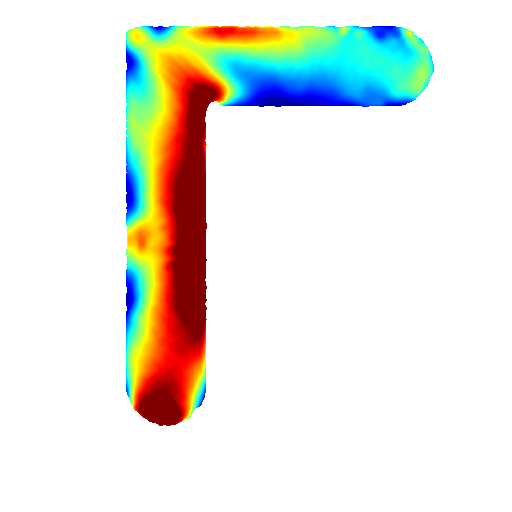}}
\put(5.5,2){\includegraphics[width=\unitlength]{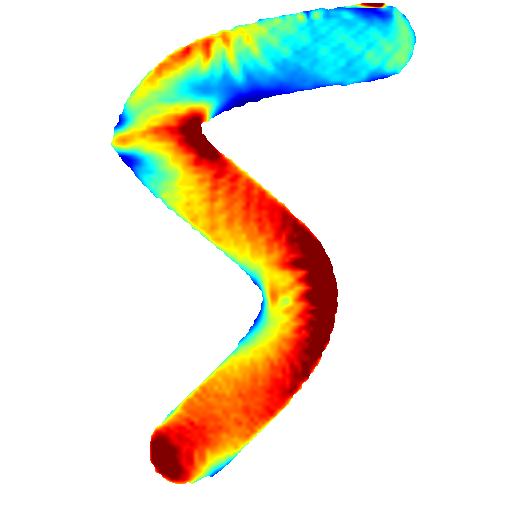}}
\put(9.5,2){\includegraphics[width=\unitlength]{GammaSCol_S}}
\end{picture}
\end{center}
\caption{
Discrete geodesics with three, five, and nine shapes (top to bottom).
Despite large shape variations and strong local rotations,
the choice of the nonlinear energy $\W$ yields good results already for coarse time steps.}
\label{fig:GammaS}
\end{figure}

To improve the efficiency of the resulting method, a cascadic finite element approach is used which proceeds from a coarse to a fine resolution of objects and deformations on a dyadic hierarchy of
meshes. Simultaneously the resolution of the discrete geodesics can be increased in time 
(\conf Figure~\ref{fig:GammaS}). In case very large nonlinear deformations occur during the optimization, from time to time the reference objects $(\hat \object_0, \ldots, \hat \object_K)$ are replaced by the current object approximations $(\object_0, \ldots, \object_K)$ and the deformations $\phi_k$ are reset to the identity deformation $\id$ for $k=0,\ldots, K$.

\subsection{Solving the constrained optimization problem associated with $\Exp{2}$}
\label{sec:Exp2numerics}
To find the zero of \eqref{eqn:fixedPointExp2} we employ the optimality condition associated with the variational definition of $\Log{2}$. Indeed, for given $\zeta$ we introduce two deformations $\psi_1,\psi_2$ and corresponding objects $\tilde \object_2 := \left(\psi_2 \circ (\id + \zeta)\right)(\object),
\tilde \object_1 := \psi_1(\object)$. Now, we associate to the discrete path $(\object, \tilde \object_1,\tilde \object_2)$ with underlying deformations $\psi_1,\psi_2$ the energy
\begin{eqnarray*}
\mathcal{E}[\psi_1,\psi_2] &:=& \int_{\object}W(\nabla \psi_1)\,\d x +
\int_{\psi_1(\object)}W\left(\nabla(\psi_2\circ(\id+\zeta)\circ \psi_1^{-1})\right) \,\d x\,.
\end{eqnarray*}
We obtain
$\partial_{\psi_1}\mathcal{E}[\psi_1,\psi_2] =0$ as the necessary condition to ensure that
$\mathcal{E}[\psi_1,\psi_2]$ actually represents the path energy in
\eqref{eqn:discretePathEnergy} connecting the two objects $\object$ and $\object_2$ 
via a discrete geodesic of order $2$ with the intermediate object $\psi_1(\object)$.
In the notion of Section \ref{sec:exp} the property of $(\object, \tilde \object_1,\tilde \object_2)$ to be 
a discrete geodesic can be phrased as  $\psi_1(\object)=\tilde \object_1[\tilde \object_2]$.
Thus, a necessary condition for \eqref{eqn:altvariationalProblemExp2} to hold is given by 
the condition 
\begin{equation}\label{eq:EXPDiscrete}
\partial_{\psi_1}\mathcal{E}[\psi_1,\psi_2]\,|_{\psi_1=\id + \zeta} =0
\end{equation}
for the remaining unknown $\psi_2$. 
With respect to the algorithmic realization we reformulate and regularized the energy as described in Section\,\ref {sec:parageodesic} above to yield 
\begin{equation*}
\mathcal{E}^{\delta_1}[\psi_1,\psi_2]
= \int_D\chi^{\delta_1}_\object \left( W(\nabla\psi_1)
+W(\nabla(\psi_2\circ(\id+\zeta))(\nabla\psi_1)^{-1})\det\nabla\psi_1\right)\,\d x\,.
\end{equation*}
Now, the discrete counterpart of  \eqref{eq:EXPDiscrete} is the condition
\begin{align}
0&=\left.\partial_{\psi_1}\mathcal{E}^{\delta_1}[\psi_1,\psi_2]\,\right|_{\psi_1=\id + \zeta}(\theta)\nonumber\\
&= \int_D \chi^{\delta_1}_\object \Bigl( \D W(\nabla(\id+\zeta)):\nabla\theta \nonumber \\
& \qquad \quad \; -\, \D W(\nabla\psi_2\circ(\id+\zeta)): (\nabla\psi_2\circ(\id+\zeta))\, \nabla\theta \,(\Id + \nabla \zeta)^{-1}\det(\Id + \nabla \zeta)
\nonumber \\
& \qquad \quad  \; +\, W(\nabla\psi_2\circ(\id+\zeta))\, \det (\Id + \nabla \zeta) \, \tr ((\Id + \nabla \zeta)^{-1} \nabla \theta)
\Bigr)\,\d x\,, \label{eq:Edelta}
\end{align}
where $A:B=\tr(A^TB)$ for matrices $A,B\in\R^{2,2}$.
This equation has to hold for all test deformations $\theta$.
In our finite element context, the corresponding test functions are taken to be all finite element basis functions so that 
\eqref{eq:Edelta} becomes a system of nonlinear equations 
which is solved for $\psi_2$ via Newton's method.
Here too, we first find $\psi_2$ on a coarse grid and then use the result as initialization of Newton's method on finer grids.

\section{Conclusions and outlook}
\label{sec:conclusions}
Based on a variational time discretization of geodesic paths in shape space we have proposed a novel time discrete geodesic calculus, which consists of discrete logarithmic and exponential maps, discrete parallel transport and a discrete connection. We demonstrate how to use this discrete calculus as a robust and efficient toolbox for shape morphing, shape extrapolation, and transport of shape features along paths of shapes. Although in this expository article we restricted ourselves to two-dimensional objects, the approach can be carried over to 3D viscous-fluid shapes. The concept can also be adapted to deformations of hypersurfaces and corresponding deformation energies, which measure tangential as well as normal bending stresses \cite{FrJaMo03,DrLiRuSc05}. For example, a generalization to the space of planar elastic curves and thin shell surfaces is feaible. Furthermore, instead of a metric structure induced by the viscous flow paradigm, the Wasserstein distance of optimal transport can be considered \cite{Villani2003}. It this case the time discretization of the Monge--Kantorovich problem proposed by Benamou and Brenier \cite{BeBr00}
is a possible starting point. Beyond these future directions of generalization, a theoretical foundation 
has to be established with existence and regularity results for the 
above-mentioned infinite dimensional shape spaces. Furthermore, the limit behaviour of the discrete geodesic calculus for vanishing time step size and the convergence of the discrete path energy to the corresponding
continuous path energy in the sense of  $\Gamma$-convergence has to be investigated (\conf the work by M\"uller and Ortiz on $\Gamma$-convergence of a time discrete action functional in the case of Hamiltonian systems
\cite{MuOr04}). Finally, given the notion of a time discrete transport, the relation of the curvature tensor to the parallel transport along the edges of a quadrilateral (\conf Proposition 1.5.8. in \cite{Klingenberg95}) can be used to define a time discrete curvature tensor, which then allows an exploration of the local geometry of shape space.

\bigskip 

{\bf Acknowledgment	} \hspace{2ex} Benedikt Wirth was supported by the German Science Foundation via the Hausdorff Center of Mathematics and by the Federal Ministry of Education and Research via CROP.SENSe.net.

{\small
\bibliographystyle{siam}
\bibliography{all,library,own}
}


\end{document}